\documentclass[a4paper,9pt]{article}
\usepackage{graphicx}
\usepackage{amsfonts}
\usepackage{amsmath}
\usepackage{amsthm}
\usepackage{amssymb}
\usepackage{fancyhdr}
\usepackage{mathtools}
\usepackage{indentfirst}
\usepackage{algorithm}
\usepackage[noend]{algpseudocode}
\usepackage{placeins}
\usepackage{listings}
\usepackage{lipsum}
\usepackage{caption}
\usepackage{subfigure}
\usepackage{xcolor}
\usepackage{comment}
\usepackage{url}
\usepackage{soul}
\usepackage{bm}

\def\longrightharpoonup{\relbar\joinrel\rightharpoonup}
\def\longleftharpoondown{\leftharpoondown\joinrel\relbar}

\def\longrightleftharpoons{
  \mathop{
    \vcenter{
       \hbox{
       \ooalign{
          \raise1pt\hbox{$\longrightharpoonup\joinrel$}\crcr
  	  \lower1pt\hbox{$\longleftharpoondown\joinrel$}
	}
      }
    }
  }
}

\allowdisplaybreaks

\newtheorem{theorem}{Theorem}
\newtheorem{remark}{Remark}
\newtheorem{lemma}{Lemma}
\newtheorem{prop}{Proposition}
\newtheorem{assump}{Assumption}
\newtheorem{corollary}{Corollary}

\definecolor{cream}{RGB}{222,217,201}

\providecommand{\keywords}[1]{\textbf{\small Keywords } #1}

\begin{document}
\title{{Pathwise estimates for effective dynamics: the case of nonlinear
vectorial reaction coordinates}}

\date{}
\author{
  Tony Leli{\`e}vre\,\textsuperscript{1} \and Wei Zhang\,\textsuperscript{2}
}

\footnotetext{$^1$\,Ecole des Ponts and INRIA, 77455 Marne-la-Vall{\'e}e Cedex 2, France}
\footnotetext{$^2$\,Zuse Institute Berlin, Takustrasse 7, 14195 Berlin, Germany}
\footnotetext{\,E-mail addresses: tony.lelievre@enpc.fr (T.~Leli{\`e}vre), wei.zhang@fu-berlin.de (W.~Zhang)}

\maketitle
\begin{abstract}
  Effective dynamics using conditional expectation was proposed in~\cite{effective_dynamics} to
  approximate the essential dynamics of high-dimensional diffusion processes along a
  given reaction coordinate. The approximation error of the effective dynamics
  when it is used to approximate the behavior of the original dynamics has been considered in recent years.
  As a continuation of the previous work~\cite{LEGOLL2017-pathwise},
  in this paper we obtain pathwise estimates
  for effective dynamics when the reaction coordinate function is either
  nonlinear or vector-valued.
\end{abstract}
\section{Introduction}
\label{sec-intro}
The evolution of many physical systems in biological
molecular dynamics and material science can be often modelled by diffusion
processes. The latter is a well-established mathematical model which allows us to
rigorously study the dynamical behavior of many real-world complex systems. 
Assuming the system is in equilibrium, one often refers to the reversible diffusion
process $x(s) \in \mathbb{R}^n$, which satisfies the stochastic differential equation (SDE)
  \begin{align}
    dx(s) = -a(x(s))\nabla V(x(s))\, ds + \frac{1}{\beta} (\nabla \cdot a)(x(s))\, ds +
    \sqrt{2\beta^{-1}} \sigma(x(s))\, dw(s)\,, \quad s \ge 0\,,
    \label{dynamics-a}
  \end{align}
  where $\beta > 0$, $w(s) \in \mathbb{R}^{n'}$ is
  an $n'$-dimensional Brownian motion with $n' \ge n$,
  and  both the potential $V : \mathbb{R}^n \rightarrow
  \mathbb{R}$ and the coefficient matrix $\sigma : \mathbb{R}^n
\rightarrow \mathbb{R}^{n \times n'}$ are smooth
functions. The symmetric matrix $a$ is related to $\sigma$ by $a=\sigma \sigma^T$ 
and in this work we always assume that $a$ is uniformly positive definite,
i.e.,
\begin{align}
  \sum_{1 \le i,j\le n} a_{ij}(x) \eta_i\eta_j  \ge c_1 |\eta|^2 \,,\quad \forall x \in \mathbb{R}^n\,, \, \eta \in
  \mathbb{R}^n\,,
  \label{a-elliptic}
\end{align}
for some constant $c_1 > 0$ and $|\cdot|$ denotes the usual Euclidean norm of
vectors. 
The notation $\nabla \cdot a$ denotes the $n$-dimensional vector whose 
components are 
$(\nabla\cdot a)_{i} = \sum\limits_{j=1}^n \frac{\partial a_{ij}}{\partial x_j}$, for $1
\le i \le n$.
Under mild conditions on the potential $V$, it is well known~\cite{Mattingly2002} that dynamics
(\ref{dynamics-a}) is ergodic with a unique invariant measure $\mu$, whose
probability density is given by
\begin{align}
  \frac{d\mu}{dx} = \frac{1}{Z} e^{-\beta V}\,, 
  \label{mu-density}
\end{align}
   where $Z = \int_{\mathbb{R}^n} e^{-\beta V} dx$
   denotes the normalization constant. 

In view of real applications, one often encounters the situation where on the one
hand the system is high-dimensional, i.e., $n \gg 1$, and on the other hand the essential 
behavior of the system can be characterized in a space whose dimension is much lower than $n$.
To study the behavior of the system in this scenario, one often assumes that
there is a (reaction coordinate) function $\xi : \mathbb{R}^n \rightarrow
\mathbb{R}^m$, where 
\begin{align}
\xi = (\xi_1, \xi_2, \cdots, \xi_m)^T\,, 
\label{xi}
\end{align}
for some $1 \le m < n$,
such that the essential dynamics of $x(s)$ can be captured by $\xi(x(s))$.
Various coarse-graining or model-reduction techniques have been developed in order to study the behavior of the dynamics along the reaction coordinate $\xi$.
We refer to~\cite{peptide_cmd,revisit-ftstring,Maragliano2006,kevrekidis2003} for related work in the study of molecular dynamics.

Notice that, applying Ito's formula, we immediately know that
$\xi(x(s))$ satisfies the SDE
\begin{align}
  d\xi(x(s)) =&\, (\mathcal{L}\xi)(x(s))\,ds +
      \sqrt{2\beta^{-1}} (\nabla \xi \sigma)\big(x(s)\big)\, dw(s)\,,
      \label{xi-x-sde}
\end{align}
where $\mathcal{L}$ is the infinitesimal generator of (\ref{dynamics-a}) and
$\nabla\xi$ denotes the $m \times n$ matrix whose entries are
$(\nabla\xi)_{ij} = \frac{\partial\xi_i}{\partial x_j}$, for $1\le i \le m$,
$1 \le j \le n$. However, (\ref{xi-x-sde}) is of limited use in practice, due to the fact that it still
depends on the original dynamics $x(s)$. In another word,
(\ref{xi-x-sde}) is not in a closed form and does not correspond to a
(Markovian) diffusion process in $\mathbb{R}^m$.
Given a reaction coordinate function $\xi$, the search of a coarse-grained
diffusion process in $\mathbb{R}^m$ in order to approximate $\xi(x(s))$ has been
studied in the past work~\cite{mimick_sde,effective_dynamics}. In particular,
the authors in~\cite{effective_dynamics} proposed an effective dynamics by
replacing the coefficients on the right hand side of (\ref{xi-x-sde}) by their
conditional expectations. In the following, we introduce several quantities in order to
explain the conditional expectation suggested in \cite{effective_dynamics}.

Given $z \in \mathbb{R}^m$, we define the level set
\begin{align}
  \Sigma_{z} = \Big\{\,x\in \mathbb{R}^n ~\Big|~\xi(x) = z\Big\}\,. 
\label{set-sigma-z}
\end{align}
Assuming it is nonempty, under certain conditions (see Remark~\ref{rmk-regular-xi}
in Subsection~\ref{subsec-assump}), $\Sigma_z$ is an $(n-m)$-dimensional submanifold of $\mathbb{R}^n$. 
We denote by $\nu_z$ the surface measure of the submanifold $\Sigma_z$ which is induced from the
Euclidean metric on $\mathbb{R}^n$. The probability measure $\mu_z$ on
$\Sigma_z$, which is defined by
\begin{align}
  d\mu_z(x) = \frac{1}{Q(z)} \frac{e^{-\beta V(x)}}{Z} \Big[\det(\nabla \xi\nabla
  \xi^T)(x)\Big]^{-\frac{1}{2}}\, d\nu_z(x)\,,
  \label{measure-mu-z}
\end{align}
has been studied in the previous
work~\cite{projection_diffusion,effective_dynamics,effective_dyn_2017,zhang2017} and 
will play an important role in the current work. 
In (\ref{measure-mu-z}), $Q(z)$ is given by 
  \begin{align}
    \begin{split}
      Q(z) =& \frac{1}{Z}\int_{\Sigma_z} e^{-\beta V(x)} \Big[\det(\nabla \xi\nabla
    \xi^T)(x)\Big]^{-\frac{1}{2}}\, d\nu_z(x) \\
    =& \frac{1}{Z} \int_{\mathbb{R}^n} \delta(\xi(x) - z)\,e^{-\beta
    V(x)}\, dx
    \end{split}
  \end{align}
and serves as the normalization constant.
Clearly, we have $\int_{\mathbb{R}^m} Q(z)\,dz = 1$.

With the above preparations, we can introduce the effective dynamics proposed in~\cite{effective_dynamics}. 
Specifically, we consider the dynamics $z(s) \in \mathbb{R}^m$ which satisfies the SDE
\begin{align}
  dz(s) = \widetilde{b}(z(s))\, ds + \sqrt{2\beta^{-1}}
  \widetilde{\sigma}(z(s))\,
  d\widetilde{w}(s)\,, \quad s \ge 0 \,,
  \label{eff-dynamics}
\end{align}
where $\widetilde{w}(s)$ is an $m$-dimensional Brownian motion,
and the coefficients are given by
\begin{align}
  \begin{split}
    \widetilde{b}_l(z) =& 
    \int_{\Sigma_z} 
    (\mathcal{L} \xi_l)(x)\, d\mu_z(x) = \mathbf{E}_\mu\bigg[ (\mathcal{L}\xi_l)(x)\,\Big|\,\xi(x) = z\bigg] \,,
\quad  1 \le l \le m \,, \\
\widetilde{\sigma}(z) =& \bigg[\int_{\Sigma_z} \big(\nabla\xi a\nabla\xi^T\big)(x)
    \,d\mu_z(x)\bigg]^{\frac{1}{2}} \,,
  \end{split}
\label{b-sigma-eff}
\end{align}
for $z \in \mathbb{R}^m$.
We recall that, given a positive definite symmetric matrix $X$,
$X^{\frac{1}{2}}$ denotes the unique positive definite symmetric matrix such that $X=
X^{\frac{1}{2}}X^{\frac{1}{2}}$. And
$\mathbf{E}_\mu\big[\cdot\,\big|\,\xi(x) = z\big]$ in (\ref{b-sigma-eff}) denotes the
conditional expectation with respect to the probability measure $\mu$. 
Furthermore, from \cite{effective_dynamics,effective_dyn_2017} we know 
that the effective dynamics (\ref{eff-dynamics}) is again both reversible and
ergodic with respect to the unique invariant measure $\widetilde{\mu}$ on
$\mathbb{R}^m$, whose probability density is $Q(z)$, i.e.,
\begin{align}
  d\widetilde{\mu}(z) = Q(z)\,dz\,.
\end{align}

With the effective dynamics (\ref{eff-dynamics}) at hand, 
it is natural to ask how good the SDE (\ref{eff-dynamics}) is when $z(s)$ is used to
approximate the process $\xi(x(s))$. In literature, the approximation error of the
effective dynamics has been studied using different criteria, such as
entropy decay rate~\cite{effective_dynamics,Sharma2017}, approximation of
eigenvalues~\cite{effective_dyn_2017}, and pathwise estimates~\cite{LEGOLL2017-pathwise}.  
As a continuation of the work~\cite{LEGOLL2017-pathwise},
in the current paper we study pathwise estimates of the effective dynamics.
While we are interested in the general case when the function $\xi$ is
nonlinear, we mention three concrete examples when  
the function
\begin{align}
  \xi(x) = (x_1, x_2, \cdots, x_m)^T\,,\quad \forall\, x \in \mathbb{R}^n
  \label{xi-linear}
\end{align}
is a linear map, since they provide useful insights and strongly motivate our current study.   
For this purpose, let us denote by $I_{m\times m}$ the identity matrix of size $m$ and write the
state $x \in \mathbb{R}^n$ as $x=(z, y) \in \mathbb{R}^{m} \times
\mathbb{R}^{n-m}$ where $y=(y_{m+1}, y_{m+2}, \cdots,
y_n)^T$, i.e., the components are indexed from $m+1$ to $n$.
Also let $\epsilon, \delta$ denote two small parameters such that $0 < \epsilon, \delta \ll 1$.
The following three cases are of particular interest.
\begin{enumerate}
  \item
    Matrix $\sigma=I_{n\times n}$ and $V(z,y) = V_0(z,y) +\frac{1}{\epsilon} V_1(y)$, where $V_0, V_1$  are two
    potential functions and $0 < \epsilon \ll 1$. 
    The SDE (\ref{dynamics-a}) becomes 
  \begin{align}
    \begin{split}
    dz_i(s) =& -\frac{\partial V_0}{\partial z_i}\big(z(s), y(s)\big)\,ds + \sqrt{2\beta^{-1}} dw_i(s)\,,\quad 1 \le i \le m\,, \\
    dy_j(s) =& -\frac{\partial V_0}{\partial y_j}\big(z(s), y(s)\big)\,ds
      -\frac{1}{\epsilon}\frac{\partial V_1}{\partial y_j}\big(y(s)\big)\,ds + \sqrt{2\beta^{-1}}\,dw_j(s)\,,\quad m+1 \le j \le n \,.
    \end{split}
    \label{linear-case-1}
  \end{align}
  \item
    Matrix $\sigma$ is constant and is given by 
  \begin{align}
  \sigma \equiv \begin{pmatrix}
    I_{m\times m} & 0 \\
    0 & \frac{1}{\sqrt{\delta}} I_{(n-m)\times (n-m)}
  \end{pmatrix}
  \,.
  \label{intro-mat-a}
\end{align}
 The SDE (\ref{dynamics-a}) becomes 
  \begin{align}
    \begin{split}
    dz_i(s) =& -\frac{\partial V}{\partial z_i}\big(z(s), y(s)\big)\,ds + \sqrt{2\beta^{-1}}
    dw_i(s)\,,\quad 1 \le i \le m\,, \\
    dy_j(s) =& -\frac{1}{\delta}\frac{\partial V}{\partial y_j}\big(z(s), y(s)\big)\,ds +
    \sqrt{\frac{2\beta^{-1}}{\delta}}\,dw_j(s)\,,\quad m+1 \le j \le n \,.
  \end{split}
  \label{linear-case-2}
  \end{align}
\item Matrix $\sigma$ is given in (\ref{intro-mat-a}) and $V(z,y) = V_0(z,y) +\frac{1}{\epsilon} V_1(y)$.
The SDE (\ref{dynamics-a}) becomes 
  \begin{align}
    \begin{split}
      dz_i(s) =& -\frac{\partial V_0}{\partial z_i}\big(z(s), y(s)\big)\,ds + \sqrt{2\beta^{-1}}
    dw_i(s)\,,\quad 1 \le i \le m\,, \\
      dy_j(s) =& -\frac{1}{\delta}\frac{\partial V_0}{\partial y_j}\big(z(s), y(s)\big)\,ds
    - \frac{1}{\epsilon\delta}\frac{\partial V_1}{\partial y_j}\big(y(s)\big)\,ds +
    \sqrt{\frac{2\beta^{-1}}{\delta}}\,dw_j(s)\,,\quad m+1 \le j \le n \,.
    \end{split}
    \label{linear-case-3}
  \end{align}
\end{enumerate}
  Among the above three cases, dynamics (\ref{linear-case-2}) in the
  second case is probably familiar, since it belongs to the typical slow-fast dynamics that has been widely studied 
  using the standard averaging technique~\cite{pavliotis2008_multiscale,liu2010,vanden-eijnden2003}.
  In this case, the probability measure $\mu_z$ in (\ref{measure-mu-z}) 
  is simply the invariant measure of the fast process $y(s)$ in
  (\ref{linear-case-2}) when $z(s)=z$ is
  fixed. Denoting by $\mathcal{L}_0$
  the infinitesimal generator of the fast process in (\ref{linear-case-2}), we
  emphasize that the decomposition of the infinitesimal generator $\mathcal{L}$ as 
  \begin{align}
  \mathcal{L} = \mathcal{L}_0 + \mathcal{L}_1
  \label{decomp-l}
\end{align} 
 plays an important role in order to derive convergence results of the system (\ref{linear-case-2})
when $\delta \rightarrow 0$~\cite{pavliotis2008_multiscale}. 
  
  With the above observation on the concrete examples in mind, let us discuss
  three key ingredients
  of our approach, which enable us to obtain pathwise estimates of the
  effective dynamics for a general reaction coordinate function $\xi$, and in particular to provide a uniform treatment of the above three examples.
  Firstly, in analogy to the averaging technique in the study of SDE
  (\ref{linear-case-2}), given the SDE (\ref{dynamics-a}) and a nonlinear
  vectorial function $\xi$, 
  we will make use of a similar decomposition of $\mathcal{L}$ to
  (\ref{decomp-l}), such that $\mathcal{L}_0$ corresponds to a diffusion process on $\Sigma_z$ whose invariant
  measure is $\mu_z$, for all $z \in \mathbb{R}^m$. 
  Secondly, for each $z \in \mathbb{R}^m$, we introduce the Dirichlet form
  $\mathcal{E}_z$ corresponding to
  $\mathcal{L}_0$ and $\mu_z$ on $\Sigma_z$. The separation of time scales in the
  dynamics (\ref{dynamics-a}) will be quantified by Poincar{\'e} inequality
  of the Dirichlet form $\mathcal{E}_z$. Thirdly, 
as the function $\xi$ is assumed to be vector-valued,
 we apply Lieb's concavity theorem~\cite{LIEB1973,ANDO1979} for positive
 definite symmetric matrices in order to get an estimate of the difference
 between two matrices in the noise term of SDEs.
  Combining these three ingredients together, we are able to generalize the pathwise estimates 
  of~\cite{LEGOLL2017-pathwise} to the case when the reaction coordinate
  function $\xi$ is either nonlinear or vector-valued. Since the time scale
  separation in the above three slow-fast examples can be characterized by the
  Poincar{\'e} inequality of $\mathcal{E}_z$ in a uniform way, our pathwise estimate results can be applied to SDEs (\ref{linear-case-1}),(\ref{linear-case-2}), and (\ref{linear-case-3}).

The rest of the paper is organized as follows. In Section~\ref{sec-main-results}, 
after introducing necessary notations and assumptions, we state our main
pathwise estimate results of this paper.
In Section~\ref{sec-time-scale}, we apply our pathwise
estimate results to three different cases when there is
a separation of time scales in the system (\ref{dynamics-a}). 
These cases are generalizations of the examples (\ref{linear-case-1}),
(\ref{linear-case-2}), and (\ref{linear-case-3}), respectively.
In Section~\ref{sec-preliminary-estimate}, we 
prove a preliminary pathwise estimate result.
Section~\ref{sec-pathwise-approx} is devoted to the proof of pathwise
estimates of effective dynamics, following the approach developed in~\cite{LEGOLL2017-pathwise}. 
In Appendix~\ref{sec-app-1}, we consider the situation when there is a coordinate
transformation such that the nonlinear reaction coordinate function can be locally reduced to a linear one. 
Appendix~\ref{sec-app-2} contains the proof of
Lemma~\ref{lemma-spectral-gap-stiff-v} in Section~\ref{sec-time-scale}. 
Finally, an error estimate of marginals under dissipative assumption is included in Appendix~\ref{sec-app-3}.

\section{Notations, assumptions, and main results}
\label{sec-main-results}
\subsection{Notations}
\label{subsec-notations}
Let us further introduce some useful notations and quantities. The infinitesimal generator of the diffusion process $(\ref{dynamics-a})$ is given by~\cite{oksendalSDE}
\begin{align}
  \begin{split}
  \mathcal{L} =& -a \nabla V \cdot \nabla + \frac{1}{\beta} (\nabla\cdot
  a) \cdot \nabla + \frac{1}{\beta} a:\nabla^2\\
  =&
\frac{e^{\beta V}}{\beta} \sum_{1 \le i, j \le n} \frac{\partial}{\partial x_i}\Big(a_{ij}
    e^{-\beta V}\frac{\partial}{\partial x_j}\Big)\,,
  \end{split}
  \label{generator-l}
\end{align}
with the notation $a:\nabla^2 = \sum\limits_{1 \le i,j \le n} a_{ij} \frac{\partial^2}{\partial
x_i\partial x_j}$.
   Given two functions $f, h : \mathbb{R}^n \rightarrow \mathbb{R}$, we define
   the weighted inner product 
  \begin{align}
    \langle f, h\rangle_\mu = \int_{\mathbb{R}^n} f(x)\,h(x)\,d\mu(x)\,,
    \label{f-h-ip}
\end{align}
whenever the right hand side exists. 
Using integration by parts, it is easy to verify that $\mathcal{L}$ 
is a self-adjoint operator with respect to the inner product (\ref{f-h-ip}). The Dirichlet form of $\mathcal{L}$ is
defined as~\cite{bakry2013analysis,villani2008optimal}
\begin{align}
 \mathcal{E}(f,h) 
 := -\langle \mathcal{L}f, h \rangle_\mu = -\langle f, \mathcal{L}h
 \rangle_\mu=\frac{1}{\beta} \int_{\mathbb{R}^n} (a\nabla f) \cdot \nabla h\, d\mu \,,
\label{lfg-flg}
\end{align}
for all functions $f, h\in \mbox{Dom}(\mathcal{L})$. 

In this work, we assume that the reaction coordinate function $\xi :
\mathbb{R}^n \rightarrow \mathbb{R}^m$ defined in (\ref{xi}) is $C^3$ smooth. 
Given $z \in \mathbb{R}^m$ and $x \in \Sigma_z$, 
we define the $m\times m$ matrix $\Phi=\nabla\xi a\nabla\xi^T$, i.e.,
\begin{align}
      \Phi_{ij} =\nabla \xi_{i} \cdot (a\nabla\xi_{j})\,,\quad \forall\,1 \le i, j \le m\,.
      \label{mat-phi-ij}
\end{align}
With a slight abuse of notation, in (\ref{mat-phi-ij}) we have denoted by $\nabla\xi_i$ the usual gradient of
the function $\xi_i$, for $1 \le i \le m$.  
Assuming the vectors $\nabla\xi_1, \nabla\xi_2, \cdots, \nabla\xi_m$ are
linearly independent (see Assumption~\ref{assump-0} in
Subsection~\ref{subsec-assump}), we have that $\Phi$ is positive definite and
therefore invertible. In this case, we denote by $A$ the positive definite symmetric matrix given by
\begin{align}
  A(x) = \big(\nabla\xi a \nabla\xi^T\big)^{\frac{1}{2}}(x) = \Phi^{\frac{1}{2}}(x)
  \,, \quad \forall~x \in \mathbb{R}^n\,,
  \label{mat-a}
\end{align}
 and we introduce the $n\times n$ matrix 
  \begin{align}
    \begin{split}
      \Pi =\,& I -  \sum_{1 \le i, j \le m} (\Phi^{-1})_{ij} \nabla\xi_{i} \otimes (a
		  \nabla\xi_{j})\,,
    \end{split}
    \label{proj-pi}
\end{align} 
where $I=I_{n\times n}$ is the identity matrix of size $n$ and $\otimes$
is the tensor product of two vectors. 
$T_x\Sigma_z$ denotes the tangent space of the submanifold
$\Sigma_z$ at $x$, and $P : T_x\mathbb{R}^n \rightarrow T_x\Sigma_z$ is the
orthogonal projection operator. It is straightforward to verify that $\Pi$ satisfies
 \begin{align}
   \begin{split}
     &\Pi^2 =\,\Pi\,, \quad \Pi^T a = a\Pi\,,\quad \Pi P=\Pi\,,\\
     & \Pi\nabla\xi_i= 0\,, \quad \Pi^T \eta =\eta \,,\quad |\Pi\eta| \ge |\eta|\,,
   \end{split}
   \label{mat-pi}
 \end{align}
 for $\forall\,\eta\, \in T_x\Sigma_z$ and $1 \le i \le m$, where the last
 inequality follows from the fact
\begin{align}
  |\Pi \eta|\,|\eta| \ge (\Pi\eta)\cdot \eta = \eta\cdot (\Pi^T\eta) =
  |\eta|^2\,, \qquad  \forall\,\eta\, \in T_x\Sigma_z\,.
\end{align}
Therefore, $\Pi^T$ defines a skew projection operator from $T_x\mathbb{R}^n$ to 
 $T_x\Sigma_z$ and it coincides with $P$ if and only if $a=I$. 

 With the matrix $\Pi$ and the expression of $\mathcal{L}$ in
 (\ref{generator-l}), we can observe that $\mathcal{L}$ can be decomposed as 
  \begin{align}
    \mathcal{L} = \mathcal{L}_0 + \mathcal{L}_1\,,
    \label{l-eqn-l0-l1}
  \end{align}
  where 
  \begin{align}
    \begin{split}
      \mathcal{L}_0 =& \frac{e^{\beta V}}{\beta}  \sum_{1\le i,j \le n} \frac{\partial}{\partial x_i} \Big(e^{-\beta V}
    (a \Pi)_{ij} \frac{\partial}{\partial x_j}\Big)\,,\\
    \mathcal{L}_1 =& \frac{e^{\beta V}}{\beta} \sum_{1\le i,j \le n}
    \frac{\partial}{\partial x_i} \Big(e^{-\beta V}
      (a(I- \Pi))_{ij} \frac{\partial}{\partial x_j}\Big)\,.
    \end{split}
    \label{l-0-1}
  \end{align}
  As already mentioned in the Introduction, an important property of the 
  decomposition (\ref{l-eqn-l0-l1})-(\ref{l-0-1}) is that, for each $z \in \mathbb{R}^m$, the operator $\mathcal{L}_0$ 
defines a diffusion process on the submanifold $\Sigma_z$ whose invariant measure
is $\mu_z$ defined in (\ref{measure-mu-z}). Furthermore, we have
\begin{align}
  \int_{\Sigma_z} \mathcal{L}_0 (fh)\, d\mu_z = 0\,,\quad \mbox{and}\quad 
  \int_{\Sigma_z} (\mathcal{L}_0f) h\, d\mu_z = \int_{\Sigma_z} f
  (\mathcal{L}_0 h)\, d\mu_z 
  \label{int-by-parts}
\end{align}
for any two smooth and bounded functions $f, h : \Sigma_z \rightarrow \mathbb{R}$.
Notice that, in the above and below, we will adopt the same notations for
both functions on $\Sigma_z$ and their smooth extensions to
$\mathbb{R}^n$.
We refer to~\cite{zhang2017} for more details. 
Corresponding to the Dirichlet form $\mathcal{E}$ in (\ref{lfg-flg}), 
we denote by $\mathcal{E}_z$ the Dirichlet form of the operator
$\mathcal{L}_0$ on $\Sigma_z$, i.e.,
\begin{align}
    \mathcal{E}_z(f,h) = -\int_{\Sigma_z} (\mathcal{L}_0 f)\, h\, d\mu_z\,, 
    \label{dirichlet-form-ez}
  \end{align}
  for all $f, h : \Sigma_z \rightarrow \mathbb{R}$ and $f, h \in
  \mbox{Dom}(\mathcal{L}_0)$. 
  Then, (\ref{mat-pi}), (\ref{l-0-1}) and (\ref{int-by-parts}) imply that 
  \begin{align}
    \mathcal{E}_z(f,h) 
    = -\int_{\Sigma_z} (\mathcal{L}_0 f) \, h\, d\mu_z\,
    = -\int_{\Sigma_z} f (\mathcal{L}_0 h) \, d\mu_z\,
    = \frac{1}{\beta} \int_{\Sigma_z} (a\Pi\nabla f) \cdot (\Pi \nabla h)\, d\mu_z\,,
    \label{dirichlet-form-z}
  \end{align}
  where the last expression is independent of the extensions $f,h$ we choose. 

On a final note, $\|X\|_{F}=\sqrt{\mbox{tr}(X^TX)}$ denotes the Frobenius norm of a matrix $X$.
Notations $\mathbf{E}_{\mu}$, $\mathbf{E}_{\mu_z}$,
$\mathbf{E}_{\widetilde{\mu}}$ will denote the mathematical expectations on
the spaces $\mathbb{R}^n$, $\Sigma_z$, and $\mathbb{R}^m$ with respect to
the probability measures $\mu$, $\mu_z$, and $\widetilde{\mu}$, respectively. By
contrast, $\mathbf{E}$ is reserved for the mathematical expectation of
paths of the dynamics (\ref{dynamics-a}) starting from $x(0)\sim\mu$. 

\subsection{Assumptions}
\label{subsec-assump}
The following assumptions will be used in the current work.
\vspace{0.1cm}
\begin{assump}
  The function $\xi : \mathbb{R}^n \rightarrow \mathbb{R}^m$ is onto, $C^3$
  smooth, and satisfies that $\mbox{rank}(\nabla\xi) = m$ at each $x \in \mathbb{R}^n$.
  \label{assump-0}
\end{assump}
\begin{remark}
  Using the terminology of differential manifold, the map $\xi$ satisfying 
  the condition $\mbox{rank}(\nabla\xi) = m$ at each point is
  called a submersion from $\mathbb{R}^n$ to $\mathbb{R}^m$.
  The assumption that $\xi$ maps onto $\mathbb{R}^m$ implies that the set
  $\Sigma_z$ in (\ref{set-sigma-z}) is nonempty for all $z \in \mathbb{R}^m$. 
  Furthermore, according to the regular value theorem~\cite{banyaga2004lectures}, $\Sigma_z$ is an $(n-m)$-dimensional submanifold of $\mathbb{R}^n$.
\label{rmk-regular-xi}
\end{remark}
\begin{assump}
  $\exists$ $L_b, L_\sigma > 0$, such that 
  \begin{align}
    |\widetilde{b}(z) - \widetilde{b}(z')| \le L_b|z-z'| \,, \quad 
    \big\|\widetilde{\sigma}(z) - \widetilde{\sigma}(z')\big\|_F \le
    L_\sigma|z-z'|\,,\qquad \forall z, z' \in \mathbb{R}^m\,. 
    \label{b-sigma-lip}
  \end{align}
  \label{assump-1}
\end{assump}
\begin{assump}
  For the matrix-valued functions $\Pi$ and $A$ defined in (\ref{proj-pi}) and
  (\ref{mat-a}) respectively, we have 
  \begin{align}
    \begin{split}
    \kappa_1^2 :=& \sum\limits_{i=1}^m \int_{\mathbb{R}^n} (\Pi\nabla \mathcal{L}\xi_i)\cdot (a\Pi\nabla
  \mathcal{L}\xi_i)\, d\mu < +\infty\,,\\
  \kappa_2^2 :=& \sum\limits_{1 \le i,j \le m} \int_{\mathbb{R}^n} (\Pi\nabla
  A_{ij})\cdot (a\Pi\nabla A_{ij})\, d\mu < +\infty\,.
  \end{split}
\end{align}
  \label{assump-2}
\end{assump}
\begin{assump}
For all $z \in \mathbb{R}^m$,
  the probability measure $\mu_z$ and the Dirichlet form $\mathcal{E}_z$ satisfy
    the Poincar{\'e} inequality with a uniform constant $\rho > 0$, i.e.,
       \begin{align}
\mbox{Var}_{\mu_z}(f) := \int_{\Sigma_z} f^2\,d\mu_z -
	 \Big(\int_{\Sigma_z} f\,d\mu_z\Big)^2 \le \frac{1}{\rho} \mathcal{E}_z(f,f)\,,
       \end{align}
       for all $\,f : \Sigma_z\rightarrow \mathbb{R}$ such that $\mathcal{E}_z(f,f) < +\infty$. 
  \label{assump-4}
\end{assump}
\vspace{0.1cm}

When studying the process (\ref{dynamics-a}) under fixed initial
condition, we also assume the following assumption. 
\vspace{0.1cm}

\begin{assump}
The Dirichlet form $\mathcal{E}$ satisfies the Poincar{\'e} inequality with constant $\alpha > 0$, i.e.,
\begin{align}
  \mbox{Var}_{\mu}(f) = \int_{\mathbb{R}^n} f^2\,d\mu -
  \Big(\int_{\mathbb{R}^n} f\,d\mu\Big)^2 \le \frac{1}{\alpha} \mathcal{E}(f,f)
	 \label{poincare-ineq-xs}
       \end{align}
  holds for all functions $f : \mathbb{R}^n \rightarrow \mathbb{R}$ such that $\mathcal{E}(f,f) < +\infty$. 
  \label{assump-5}
\end{assump}

\subsection{Main results}
\label{subsec-main-results}
In order to state our pathwise estimate results,
we need to first construct a version of the effective dynamics $z(s)$, such
that the Brownian motion $\widetilde{w}(s)$ in (\ref{eff-dynamics}) is coupled to the Brownian
motion $w(s)$ in the original dynamics (\ref{dynamics-a}).
For this purpose, we introduce the process $\widetilde{w}(s)$ which satisfies
\begin{align}
  d\widetilde{w}(s) = \big(A^{-1} \nabla\xi \sigma\big)(x(s))\, dw(s)=
  \Big[(\nabla\xi a \nabla\xi^T)^{-\frac{1}{2}} \nabla\xi \sigma\Big](x(s))\, dw(s)\,, \quad s \ge 0\,.
  \label{tilde-w-coupled}
\end{align}
Using L{\'e}vy's characterization of Brownian motion and the relation $a=\sigma\sigma^T$, it is straightforward to verify that
(\ref{tilde-w-coupled}) indeed defines an $m$-dimensional
Brownian motion.
With this choice of the driving noise, the effective dynamics (\ref{eff-dynamics}) becomes
\begin{align}
  \begin{split}
  dz(s) =&\, \widetilde{b}(z(s))\, ds + \sqrt{2\beta^{-1}}
  \widetilde{\sigma}(z(s))\, d\widetilde{w}(s)\\
  =&\, \widetilde{b}(z(s))\, ds + \sqrt{2\beta^{-1}}
    \widetilde{\sigma}(z(s))\Big[\big(\nabla\xi a \nabla\xi^T\big)^{-\frac{1}{2}}
  \nabla\xi\sigma\Big](x(s))\, dw(s)\,.
  \end{split}
  \label{modified-eff-dynamics}
\end{align}
Accordingly, the difference between $z(s)$ and $\xi(x(s))$ satisfies 
\begin{align}
  \begin{split}
     d\big(\xi(x(s)) - z(s)\big) = &\,
  \Big[(\mathcal{L}\xi)(x(s)) - \widetilde{b}(z(s))\Big]\,ds 
       + \sqrt{2\beta^{-1}} \Big[\big(\nabla\xi a
      \nabla\xi^T\big)^{\frac{1}{2}}(x(s)) - 
	\widetilde{\sigma}(z(s))\Big]\,d\widetilde{w}(s)\,\\
	= &\, \varphi(x(s))\,ds + \Big[\widetilde{b}\big(\xi(x(s))\big) - \widetilde{b}(z(s))\Big]\,ds 
       + \sqrt{2\beta^{-1}} \Big[A(x(s)) - \widetilde{\sigma}(z(s))\Big]\,d\widetilde{w}(s)\,,
  \end{split}
	\label{xi-z-diff}
\end{align}
where we have introduced the function $\varphi : \mathbb{R}^n \rightarrow
\mathbb{R}^m$, given by 
\begin{align}
  \varphi(x) = (\mathcal{L}\xi)(x) - \widetilde{b}(\xi(x))\,,\quad \forall~x \in \mathbb{R}^n\,.
  \label{f-lxi-b-diff}
\end{align} 

Let us first consider the case when the dynamics $x(s)$ starts from equilibrium, i.e., $x(0)\sim \mu$.
Using relatively simple argument, in Section~\ref{sec-preliminary-estimate} we
obtain our first pathwise estimate of the effective dynamics, which is stated
below.
\begin{prop}
  Suppose that Assumptions~\ref{assump-0}-\ref{assump-4} hold. 
  $x(s)$ satisfies the SDE (\ref{dynamics-a}) starting from
  $x(0)\sim \mu$, and $z(s)$ is the effective dynamics (\ref{modified-eff-dynamics}) with $z(0) =
  \xi(x(0))$.
  For all $t \ge 0$, we have 
\begin{align}
  \mathbf{E}\Big(\sup_{0 \le s \le t}\big|\xi(x(s)) - z(s)\big|^2\Big) \le 
    \frac{3\,t}{\beta\rho}\Big(\kappa_1^2\,t + \frac{32\kappa_2^2}{\beta}\Big) e^{Lt}\,,
    \label{pathwise-1-over-rho}
\end{align}
where $L = 3L_b^2+\frac{48L_\sigma^2}{\beta}+1$.
\label{prop-pathwise-1}
\end{prop}

Following the approach of~\cite{LEGOLL2017-pathwise} and applying the forward-backward martingale
method~\cite{lyons1994,komorowski2012fluctuations}, in Section~\ref{sec-pathwise-approx}, we prove the following improved pathwise estimate of the effective dynamics.

\begin{theorem}
  Suppose that Assumptions~\ref{assump-0}-\ref{assump-4} hold. 
  $x(s)$ satisfies the SDE (\ref{dynamics-a}) starting from
  $x(0)\sim \mu$, and $z(s)$ is the effective dynamics (\ref{modified-eff-dynamics}) with $z(0) =
  \xi(x(0))$. For all $t \ge 0$, we have 
\begin{align}
  \mathbf{E}\Big(\sup_{0 \le s \le t}\big|\xi(x(s)) - z(s)\big|^2\Big) \le 
    \frac{3t}{\beta\rho}\Big(\frac{27\kappa_1^2}{2\rho} + \frac{32\kappa_2^2}{\beta} \Big) e^{Lt}\,,
  \label{pathwise-bound}
\end{align}
where $L = 3L_b^2+\frac{48L_\sigma^2}{\beta}+1$.
\label{thm-pathwise-estimate}
\end{theorem}
\begin{remark}
  We make two remarks.
  \begin{enumerate}
    \item
      Theorem~\ref{thm-pathwise-estimate} replies on global Lipschitz conditions
(Assumption~\ref{assump-1}) on the coefficients of the effective dynamics.
Alternatively, in Appendix~\ref{sec-app-3} we show that the dissipative assumption~\cite{liu2010,zhws16}
      can be exploited as well, in order to obtain estimate of
      $\mathbf{E}|\xi(x(t))-z(t)|^2$, i.e., the mean square error of the
      marginals between the two processes.
\item
  The setting of~\cite{LEGOLL2017-pathwise} corresponds to the case when
      $a=I_{n\times n}$ and the function $\xi$ is linear. In this case,
      the constant $\kappa_2=0$ and the forward-backward martingale method
      indeed improves the pathwise estimate error from $\mathcal{O}(\frac{1}{\rho})$
      to $\mathcal{O}(\frac{1}{\rho^2})$. 
  However, in general cases when either $\xi$ is nonlinear or the matrix $a$ is non-identity,
  $\kappa_2$ is typically non-zero and the error bound (\ref{pathwise-bound})
  is still $\mathcal{O}(\frac{1}{\rho})$, i.e., the same as
  Proposition~\ref{prop-pathwise-1}. This is partially due to the existence of the martingale term in (\ref{xi-z-diff}). 
 Nevertheless, the dependence of the error bound (\ref{pathwise-bound}) on
 the parameter $\kappa_2$ seems necessary. And from
 Assumption~\ref{assump-2} we can observe that $\kappa_2$ will be small when the
 matrix function $A$
 in (\ref{mat-a}) is close to a constant on each submanifold $\Sigma_z$. 
  \end{enumerate}
       \label{rmk-4}
\end{remark}
\vspace{0.2cm}

Now we turn to more general initial conditions. Applying Theorem~\ref{thm-pathwise-estimate}, in
Section~\ref{sec-pathwise-approx} we will prove the following pathwise estimate result.

\begin{theorem}
Suppose that Assumptions~\ref{assump-0}-\ref{assump-4} hold. 
  $x(s)$ satisfies the SDE (\ref{dynamics-a}) and $z(s)$ is the effective
  dynamics (\ref{modified-eff-dynamics}) with $z(0) = \xi(x(0))$. 
  \begin{enumerate}
    \item
  Suppose $x(0)\sim \bar{\mu}$, where the probability measure $\bar{\mu}$ is absolutely continuous with respect to
  $\mu$ such that $$\int_{\mathbb{R}^n} \Big(\frac{d\bar{\mu}}{d\mu}\Big)^2\,
      d\mu < +\infty\,.$$ 
  Then, for all $t \ge 0$, we have 
\begin{align}
  \mathbf{E}\bigg(\sup_{0 \le s \le t}\big|\xi(x(s)) - z(s)\big|~\Big|~x(0)
  \sim \bar{\mu}\bigg) \le 
  \sqrt{t}
  \Big(
  \frac{9\kappa_1}{\sqrt{2\beta}\rho} +
  \frac{12\kappa_2}{\beta\sqrt{\rho}}\, \Big) 
  \Big[\int_{\mathbb{R}^n} \Big(\frac{d\bar{\mu}}{d\mu}\Big)^{2} d\mu\Big]^{\frac{1}{2}} 
e^{Lt}\,,
  \label{pathwise-bound-mu-bar}
\end{align}
      where $L = \frac{3}{2}L_b^2+\frac{24L_\sigma^2}{\beta}+\frac{1}{2}$.
    \item
      Suppose $x(0)=x' \in \mathbb{R}^n$ is fixed and that
      Assumption~\ref{assump-5} holds. Both the function $\varphi$ in
      (\ref{f-lxi-b-diff}) and the matrix $A$ in (\ref{mat-a}) are bounded on
      $\mathbb{R}^n$, i.e., $|\varphi(x)| \le C_1$ and $\|A(x)\|_F\le C_2$, $\forall x
      \in \mathbb{R}^n$, for some $C_1, C_2 > 0$. Let $\mu_s$ be the
      probability measure of $x(s)$ for $s\ge 0$ and denote
      $p_s=\frac{d\mu_s}{d\mu}$ when $s>0$. Given any $0 < t_0 \le t_1 \le t$, we have 
\begin{align}
  \begin{split}
    & \mathbf{E}\bigg(\sup_{0 \le s \le t}\big|\xi(x(s)) -
  z(s)\big|~\Big|~x(0)=x'\bigg) \\
	    \le & \bigg\{\sqrt{t}
  \Big(
  \frac{9\kappa_1}{\sqrt{2\beta}\rho} +
  \frac{12\kappa_2}{\beta\sqrt{\rho}}\, \Big) 
	    \bigg[1 + e^{-\alpha (t_1 - t_0)} \Big(\int_{\mathbb{R}^n} p^2_{t_0}
	    d\mu\Big)^{\frac{1}{2}}\bigg]  
	     + \sqrt{t_1}\Big(3C_1 \sqrt{t_1} +
	     \frac{18C_2}{\sqrt{\beta}}\Big)\bigg\}  e^{Lt} \,,
  \end{split}
  \label{pathwise-bound-fix-x0}
\end{align}
	where $L = \frac{3}{2}L_b^2+\frac{24L_\sigma^2}{\beta}+\frac{1}{2}$
	and $\alpha$
      is the Poincar{\'e} constant in (\ref{poincare-ineq-xs}).
\end{enumerate}
\label{thm-pathwise-estimate-mu-bar}
\end{theorem}

\begin{remark}
Notice that $\mu_0$ in Theorem~\ref{thm-pathwise-estimate-mu-bar} will be a
delta measure when $x(s)$ starts from a fixed initial condition $x(0)=x'$. 
The time $t_0>0$ is introduced to make sure that $\int_{\mathbb{R}^n} p_{t_0}^2 d\mu<+\infty$. 
We point out that  Assumption~\ref{assump-5} is not needed when $t_1 = t_0$, and
in this case (\ref{pathwise-bound-fix-x0}) becomes 
\begin{align}
  \begin{split}
    & \mathbf{E}\bigg(\sup_{0 \le s \le t}\big|\xi(x(s)) -
  z(s)\big|~\Big|~x(0)=x'\bigg) \\
	    \le & \bigg\{\sqrt{t}
  \Big(
  \frac{9\kappa_1}{\sqrt{2\beta}\rho} +
  \frac{12\kappa_2}{\beta\sqrt{\rho}}\, \Big) 
	    \bigg[1 + \Big(\int_{\mathbb{R}^n} p^2_{t_0}
	    d\mu\Big)^{\frac{1}{2}}\bigg]  
	     + \sqrt{t_0}\Big(3C_1 \sqrt{t_0} +
	     \frac{18C_2}{\sqrt{\beta}}\Big)\bigg\}  e^{Lt} \,.
  \end{split}
  \label{pathwise-bound-fix-x0-t1-t2-same}
\end{align}
Comparing to (\ref{pathwise-bound-fix-x0-t1-t2-same}), 
the estimate (\ref{pathwise-bound-fix-x0}) allows us to 
  further optimize the upper bound of the error estimate by varying $t_1\in [t_0,t]$,
under Assumption~\ref{assump-5}. 
       \label{rmk-5}
\end{remark}

In the next section, we will apply Theorem~\ref{thm-pathwise-estimate} to 
three different scenarios when there is time scale separation in the system.
We refer to Corollary~\ref{corollary-case-1}-\ref{corollary-case-3} in Section~\ref{sec-time-scale} for the pathwise estimate result in each case.  
\section{Separation of time scales in diffusion processes}
\label{sec-time-scale}
Our pathwise estimates of the effective dynamics rely on Assumption 
\ref{assump-4}, which characterizes the existence of the time scale separation in the system (\ref{dynamics-a}). 
In this section, we consider the relation between 
the structure of the SDE (\ref{dynamics-a}) and 
the emergence of the time scale separation phenomena in the process $x(x)$.
We apply our pathwise estimates to different scenarios, and in particular we obtain pathwise estimates of the
effective dynamics for the SDEs (\ref{linear-case-1}), (\ref{linear-case-2}), and (\ref{linear-case-3}) in the Introduction.

  Roughly speaking, the time scales of the dynamics (\ref{dynamics-a}) are related to 
  the magnitudes of coefficients in the infinitesimal generator
  $\mathcal{L}$.
  With the choice of the reaction coordinate function $\xi$ in (\ref{xi}) and the corresponding
  decomposition (\ref{l-eqn-l0-l1})-(\ref{l-0-1}) of the infinitesimal generator $\mathcal{L}$, 
  we are interested in cases when 
  \begin{align}
    \textit{operator $\mathcal{L}_0$ contains large coefficients, while $\mathcal{L}_1$ does not.}
    \hspace{0.5cm} \label{scale-separation} 
  \end{align}

  \noindent As we will see,  
  condition (\ref{scale-separation}) often implies that the operator $\mathcal{L}_0$ has a large
  spectral gap while the process $\xi(x(s))$ evolves relatively slowly.
  From the expression of $\mathcal{L}_0$ in (\ref{l-0-1}), we can observe that large
  coefficients in $\mathcal{L}_0$ may come from either the potential $V$
  or the (eigenvalues of) matrix $a$.  Motivated by the concrete examples (\ref{linear-case-1}), (\ref{linear-case-2}), and (\ref{linear-case-3}) in the Introduction, 
in the following we consider three different cases. For simplicity, we will assume the existence of small
  parameters $\epsilon$ or $\delta$ whose specific values are not necessarily known,
  such that the magnitudes of small and large quantities
  correspond to $\mathcal{O}(1)$ and $\mathcal{O}(\frac{1}{\epsilon})$ (or
  $\mathcal{O}(\frac{1}{\delta})$),
  respectively.
  \vspace{0.1cm}

  \textbf{Case $1$}. In the first case, let us assume that the potential $V$
  contains stiff components of $\mathcal{O}(\frac{1}{\epsilon})$, while the eigenvalues of matrix $a$ are $\mathcal{O}(1)$. 
  From expressions in (\ref{l-0-1}), we see that the condition (\ref{scale-separation}) holds if 
     \begin{align}
       \Pi^Ta\nabla V = \mathcal{O}\big(\frac{1}{\epsilon}\big)\,,\quad
       \mbox{and}\quad
       (I-\Pi^T)a\nabla V= \mathcal{O}(1)\,. 
       \label{pi-a-grad-v-large}
     \end{align}
     Since $\Pi^T$ is a skew projection operator from $T_x\mathbb{R}^n$ to
     $T_x\Sigma_z$ at each $x\in \mathbb{R}^n$, (\ref{pi-a-grad-v-large}) is
     equivalent to that the stiff component of $a\nabla V$ lies in the
     subspace $T_x\Sigma_z$ at each $x$. 
     As a concrete example, assume that the potential $V$ takes the form 
  \begin{align}
  V(x) = V_0(x) +\frac{1}{\epsilon} V_1(x)\,,
  \label{v-v0-v1}
  \end{align}
  where $V_0, V_1 : \mathbb{R}^n \rightarrow \mathbb{R}$ are two potential
  functions, $\epsilon > 0$ is a small parameter and that the condition
     \begin{align}
     (I-\Pi)^T a\nabla
   V_1\equiv 0 
    \label{decomp-pi-orthgonal-v1}
 \end{align}
 is satisfied at each $x$. 
     Clearly, in this case we have 
  \begin{align}
    \begin{split}
    \Pi^Ta\nabla V =&\, \Pi^Ta\nabla V_0 + \frac{1}{\epsilon} a\nabla V_1\, \\
    (I-\Pi)^Ta\nabla V =&\, (I-\Pi)^Ta\nabla V_0 \,,
    \end{split}
    \label{v-decomp-eps}
  \end{align}
     which implies that the condition (\ref{scale-separation}) holds. 
     In fact, corresponding to the potential $V$ in (\ref{v-v0-v1}), the
     probability measure $\mu_z$ in (\ref{measure-mu-z}) becomes 
  \begin{align}
    d\mu_z(x) = \frac{1}{Z\,Q(z)}\exp\Big[-\beta\Big(V_0(x) + \frac{1}{\epsilon}
    V_1(x)\Big)\Big] \Big[\mbox{det}(\nabla \xi\nabla
    \xi^T)(x)\Big]^{-\frac{1}{2}}\,d\nu_z(x)\,,
  \end{align}
  for each $z \in \mathbb{R}^m$.
  This measure indeed satisfies a Poincar{\'e} inequality with a large spectral
     gap if the potential $V_1$ is convex on $\Sigma_z$. Precisely, we have the following result.
     \vspace{0.1cm}

     \begin{lemma}
       Suppose the function $V_0$ in (\ref{v-v0-v1}) is bounded on $\Sigma_z$.
       $V_1$ is both $C^2$ smooth and $K$-convex on $\Sigma_z$ for some $K>0$. Matrix
       $a$ satisfies the uniform elliptic condition (\ref{a-elliptic}) with some constant $c_1 > 0$, 
        and the function $\xi$ has bounded derivatives up to order $3$. 
        Then $\exists \epsilon_0,\,C \ge 0$, which may depend on $V_0$, $a$, $\xi$ and $\beta$, such that when $\epsilon \le \epsilon_0$, 
       the Poincar{\'e} inequality 
       \begin{align}
	 \mbox{Var}_{\mu_z}(f) = \int_{\Sigma_z} f^2\,d\mu_z - 
	 \Big(\int_{\Sigma_z} f\,d\mu_z\Big)^2 \le 
	 \frac{C\epsilon}{c_1 K} \mathcal{E}_z(f,f)
	 \label{poincare-ineq-prop}
       \end{align}
holds for all functions $f : \Sigma_z\rightarrow \mathbb{R}$ which satisfy $\mathcal{E}_z(f,f) < +\infty$. 
       \label{lemma-spectral-gap-stiff-v}
     \end{lemma}

     The proof of Lemma~\ref{lemma-spectral-gap-stiff-v} can be found
     in Appendix~\ref{sec-app-2}. Applying Theorem~\ref{thm-pathwise-estimate} and
     Lemma~\ref{lemma-spectral-gap-stiff-v}, we can obtain the pathwise
     estimate of the effective dynamics in this case.
     \vspace{0.1cm}

\begin{corollary}
  Assume Assumptions~\ref{assump-0}-\ref{assump-1} hold.
     Let the potential $V$ be given in (\ref{v-v0-v1}), where $\epsilon>0$ is
     a small parameter, such that the assumptions of Lemma~\ref{lemma-spectral-gap-stiff-v} hold uniformly for $z \in \mathbb{R}^m$.
     Further suppose that the condition (\ref{decomp-pi-orthgonal-v1}) is satisfied.
  Let $x(s)$ satisfy the SDE (\ref{dynamics-a}) starting from
  $x(0)\sim \mu$, and $z(s)$ be the effective dynamics (\ref{modified-eff-dynamics}) with $z(0) =
  \xi(x(0))$.
 Define $L = 3L_b^2+\frac{48L_\sigma^2}{\beta}+1$.
    Then $\exists \epsilon_0 \ge 0$, such that when $\epsilon \le \epsilon_0$, we have 
\begin{align}
  \mathbf{E}\Big(\sup_{0 \le s \le t}\big|\xi(x(s)) - z(s)\big|^2\Big) \le 
    t\Big(\frac{C_1\epsilon}{K} + \frac{C_2\epsilon^2}{K^2}\Big) e^{Lt}\,,
  \label{pathwise-bound-case-1}
\end{align}
for some constants $C_1, C_2>0$ which are independent of $\epsilon$ and $K$.
\label{corollary-case-1}
\end{corollary}
\begin{proof}
      From the definition of $\mathcal{L}$ in (\ref{generator-l}), the
      condition (\ref{decomp-pi-orthgonal-v1}), as well as the boundedness of
      both the matrix $a$ and $\nabla \xi$, we know that 
      Assumption~\ref{assump-2} holds with constants $\kappa_1, \kappa_2$
      which are independent of $\epsilon$.
Lemma~\ref{lemma-spectral-gap-stiff-v} implies that Assumption~\ref{assump-4} is met with $\rho=
\frac{c_1 K}{C\epsilon}$.  
Therefore, the estimate (\ref{pathwise-bound-case-1}) follows by applying Theorem~\ref{thm-pathwise-estimate}.
\end{proof}

     \begin{remark}
       Given $x \in \mathbb{R}^n$, in Appendix~\ref{sec-app-1} we will study
       the condition under which there exists a function $\phi : \mathbb{R}^n \rightarrow \mathbb{R}^{n-m}$ and a coordinate
       transformation $G(x) = \big(\xi(x), \phi(x)\big),$ such that $G$ is one
       to one in a neighborhood of $x$ and that 
       \begin{align}
       \nabla\xi a \nabla\phi^T \equiv 0
       \label{assump-xi-phi-repeat}
     \end{align}
     is satisfied. See the condition~(\ref{assump-xi-phi}) in Appendix~\ref{sec-app-1}. For simplicity, let us assume
that the map $\phi$ exists globally such that $G$ is one to one from $\mathbb{R}^n$
to itself. In this case,  Assuming the potential $V$ is given in (\ref{v-v0-v1})
       with $V_1(x) = \widetilde{V}_1(\phi(x))$ for some function
       $\widetilde{V}_1 : \mathbb{R}^{n-m}\rightarrow \mathbb{R}$, then 
       (\ref{decomp-pi-orthgonal-v1}) holds because of~(\ref{assump-xi-phi-repeat}),
       and the SDE of $\bar{y}(s) = \phi(x(s))$ 
       has a large drift term which involves the small parameter $\epsilon$, while the
       SDE of $\bar{z}(s) = \xi(x(s))$ is independent of $\epsilon$.
       See (\ref{xi-phi-dynamics}) in Appendix~\ref{sec-app-1} for details.
       According to Proposition~\ref{prop-phi-sde} in Appendix~\ref{sec-app-1} and
       Lemma~\ref{lemma-spectral-gap-stiff-v} above, the invariant
       measure of the dynamics $\bar{x}(s) = G^{-1}(z, \bar{y}(s)) \in
       \Sigma_z$ with $\bar{z}(s) = z$ being fixed (see (\ref{y-fix-z})) is
       $\mu_z$ and satisfies the Poincar{\'e} inequality (\ref{poincare-ineq-prop}).
       As a concrete example, consider the linear reaction coordinate
       case when 
  \begin{align}
    \xi=(x_1, x_2, \cdots, x_m)^T\,, \quad \phi=(x_{m+1}, \cdots, x_{n})^T\,,
       \label{rmk-1-xi-phi}
     \end{align}
     and $a\equiv I_{n\times n}$, where the potential function $V_1(x)=V_1(x_{m+1},
     \cdots, x_n)$ is independent of the first $m$ components of $x$.
In this case, we have 
\begin{align}
  \Phi = I_{m\times m}\,,\quad 
  \Pi = \begin{pmatrix}
    0 & 0 \\
    0 & I_{(n-m)\times (n-m)}
  \end{pmatrix}\,,
  \label{rmk-1-linear-phi-pi}
\end{align}
and the dynamics (\ref{dynamics-a}) reduces to the SDE (\ref{linear-case-1}) in the Introduction.
Correspondingly, the constant $C_1 =0$ in (\ref{pathwise-bound-case-1}), since $\kappa_2 = 0$ in
Assumption~\ref{assump-2}. Thereore, we have the pathwise estimate 
  \begin{align}
    \mathbf{E}\Big(\sup_{0 \le s \le t}\big|\xi(x(s)) - z(s)\big|^2\Big) \le
    \frac{C_2 \epsilon^2t}{K^2} e^{Lt}\,.
  \label{pathwise-bound-case-1-linear}
\end{align}
\label{rmk-1}
     \end{remark}

  \textbf{Case $2$}. In the second case, let us assume that the potential
  function $V$ is $\mathcal{O}(1)$, but the matrix $a$ has widely spread
  eigenvalues at two different orders of magnitude. 
   Specifically, suppose that the eigenvalues $\lambda_i$ of $a$ satisfy 
that $\lambda_i=\mathcal{O}(1)$ for $1 \le i \le m$, and
$\lambda_i=\mathcal{O}(\frac{1}{\delta})$ for $m+1 \le i \le n$.
     Furthermore, the reaction coordinate function
     $\xi$ is chosen in a way such that at each state $x$ the linear subspace 
     \begin{align}
       \mbox{span}\Big\{\nabla \xi_1, \nabla \xi_2, \cdots, \nabla \xi_m\Big\}
     \end{align}
     coincides with the subspace spanned by the eigenvectors of matrix $a$
     which correspond to the (small) eigenvalues $\lambda_1, \lambda_2, \cdots, \lambda_m$.  
    Since the projection matrix $\Pi^T$ satisfies (\ref{mat-pi}), 
    in this case we have 
     \begin{align} 
       \Pi^Ta = a\Pi=\mathcal{O}\Big(\frac{1}{\delta}\Big)\,, \quad \mbox{and}\quad (I-\Pi^T)a = a(I-\Pi) = \mathcal{O}(1)\,. 
       \label{pi-a-order-delta}
     \end{align}
Therefore, from expressions in (\ref{l-0-1}) we can conclude that the condition
(\ref{scale-separation}) is satisfied. 
Notice that, different from the previous case, now the probability measure
$\mu_z$ in (\ref{measure-mu-z}) does not depend on $\delta$, while the
Dirichlet form $\mathcal{E}_z$ in (\ref{dirichlet-form-z}) does. Concerning
the Poincar{\'e} inequality, we have the following straightforward result.
\vspace{0.1cm}

     \begin{lemma}
       Given $\rho_0, \delta > 0$.  Recall that $P$ is the orthogonal projection operator from $T_x\mathbb{R}^n$ to $T_x\Sigma_z$. Assume that the probability measure $\mu_z$ satisfies 
       \begin{align}
	 \mbox{Var}_{\mu_z}(f) = \int_{\Sigma_z} f^2\,d\mu_z - 
	 \Big(\int_{\Sigma_z} f\,d\mu_z\Big)^2 \le \frac{1}{\beta\rho_0}
	 \int_{\Sigma_z} |P\nabla f|^2 d\mu_z\,,
	 \label{poincare-ineq-no-eps-case-2}
       \end{align}
        for all functions $f : \Sigma_z\rightarrow \mathbb{R}$ such that
	$\int_{\Sigma_z} |P\nabla f|^2 d\mu_z<+\infty$ (after being extended
	to a function on $\mathbb{R}^n$). Also assume the matrices $a$ and
	$\Pi$ satisfy 
	\begin{align}
	  (a\Pi \eta)\cdot \eta \ge \frac{c_2}{\delta} |\eta|^2 \,,\quad
	  \forall\,\eta \in T_x\Sigma_z\,,~\forall\, x \in \Sigma_z\,,
	  \label{a-pi-elliptic-eps}
	\end{align}
	for some $c_2 > 0$, which is independent of $\delta$. Then we have the Poincar{\'e} inequality 
       \begin{align}
	 \mbox{Var}_{\mu_z}(f) = \int_{\Sigma_z} f^2\,d\mu_z - 
	 \Big(\int_{\Sigma_z} f\,d\mu_z\Big)^2 \le \frac{\delta}{c_2 \rho_0} \mathcal{E}_z(f,f)\,,
	 \label{poincare-ineq-case-2}
       \end{align}
        for all functions $f : \Sigma_z\rightarrow \mathbb{R}$ such that $\mathcal{E}_z(f,f) < +\infty$.
       \label{lemma-spectral-gap-case-2}
     \end{lemma}
     \begin{proof}
       Notice that (\ref{mat-pi}) implies $\Pi P=\Pi$ and $\Pi^Ta\Pi=a\Pi^2=a\Pi$. Since
       $P\nabla f \in T_x\Sigma_z$, using (\ref{a-pi-elliptic-eps}) we can deduce 
       \begin{align*}
	 |P\nabla f|^2 \le \frac{\delta}{c_2} (a\Pi P\nabla f)\cdot (P\nabla
	 f) = \frac{\delta}{c_2} (\Pi^Ta\Pi\nabla f)\cdot (P\nabla f) 
= \frac{\delta}{c_2} (a\Pi \nabla f)\cdot (\Pi\nabla f)\,. 
       \end{align*}
       The conclusion (\ref{poincare-ineq-case-2}) follows by recalling the
       definition of $\mathcal{E}_z$ in (\ref{dirichlet-form-z}).
     \end{proof}
     \vspace{0.1cm}
Applying Theorem~\ref{thm-pathwise-estimate} and
Lemma~\ref{lemma-spectral-gap-case-2}, we can obtain the pathwise estimate of
the effective dynamics in this case.

\begin{corollary}
Assume Assumptions~\ref{assump-0}-\ref{assump-1} hold.
Suppose that the matrix $a$ satisfies the uniform elliptic condition
(\ref{a-elliptic}) with some constant $c_1 > 0$, and the function $\xi$ has bounded derivatives up to order $2$.
Furthermore, the matrices $a$ and $\Pi$ satisfy 
	\begin{align}
 \frac{c_2}{\delta} |\eta|^2
\le (a\Pi \eta)\cdot \eta \le \frac{c_2'}{\delta} |\eta|^2 \,,\quad
	  \forall\,\eta \in T_x\Sigma_z\,,~\forall\, x \in \Sigma_z\,,~ z \in
	  \mathbb{R}^m\,,
	  \label{a-pi-elliptic-eps-upper-lower}
	\end{align}
	for some $0 < c_2\le c_2'$, which are independent of $\delta>0$. 
  Suppose $\mu_z$ satisfies the Poincar{\'e} inequality
  (\ref{poincare-ineq-no-eps-case-2}) with the constant
  $\rho_0>0$, uniformly for $z \in \mathbb{R}^m$. 
  The matrix $A$ in (\ref{mat-a}) is bounded with bounded derivatives up to order $2$.
  Let $x(s)$ satisfy the SDE (\ref{dynamics-a}) starting from
  $x(0)\sim \mu$, and $z(s)$ be the effective dynamics (\ref{modified-eff-dynamics}) with $z(0) =
  \xi(x(0))$.  Define $L = 3L_b^2+\frac{48L_\sigma^2}{\beta}+1$. 
	  We have 
  \begin{align}
    \mathbf{E}\Big(\sup_{0 \le s \le t}\big|\xi(x(s)) - z(s)\big|^2\Big) \le t \Big(\frac{C_1}{\rho_0} + \frac{C_2\delta}{\rho_0^2}\Big) e^{Lt}\,,
  \label{pathwise-bound-case-2}
\end{align}
for some constants $C_1, C_2>0$ which are independent of $\delta$ and $\rho_0$.
\label{corollary-case-2}
\end{corollary}
\begin{proof}
      Condition (\ref{a-pi-elliptic-eps-upper-lower}) implies that 
      Assumption~\ref{assump-2} holds but the constants $\kappa_1, \kappa_2$
      may depend on $\delta$ such that $\kappa_1^2,
      \kappa_2^2 \le \frac{C}{\delta}$, for some $C>0$.
      Assumption~\ref{assump-4} is met with $\rho=\frac{c_2 \rho_0}{\delta}$
      as a consequence of Lemma~\ref{lemma-spectral-gap-case-2}.
      Therefore, the estimate~(\ref{pathwise-bound-case-2}) follows by applying Theorem~\ref{thm-pathwise-estimate}.
\end{proof}

\begin{remark}
  Consider the function $\phi$ in Appendix~\ref{sec-app-1} 
  which satisfies the condition $\nabla\xi a \nabla\phi^T \equiv 0$.
     In this case, the subspace $\mbox{span}\big\{\nabla \phi_1, \nabla \phi_2, \cdots, \nabla \phi_{n-m}\big\}$
     coincides with the subspace spanned by the eigenvectors of matrix $a$
     which correspond to the large
     eigenvalues $\lambda_{m+1}$, $\lambda_{m+2}$, $\cdots$, $\lambda_{n}$. 
     And we have 
 \begin{align}
   \begin{pmatrix}
     \nabla\xi\\
     \nabla \phi
   \end{pmatrix}
   a 
   \begin{pmatrix}
     \nabla\xi^T & \nabla \phi^T 
   \end{pmatrix}
   = \begin{pmatrix}
     \nabla\xi a\nabla\xi^T & 0 \\
     0 & \nabla\phi a\nabla\phi^T 
   \end{pmatrix}
   =  \begin{pmatrix}
     \mathcal{O}(1) & 0 \\
     0 & \mathcal{O}(\frac{1}{\delta})
   \end{pmatrix}\,.
 \end{align}
 Therefore, we can observe that the dynamics of $\phi(x(s))$  will be
 fast, while the dynamics of $\xi(x(s))$ is relatively slow. See the equation (\ref{xi-phi-dynamics})
 in Appendix~\ref{sec-app-1}.
  As a concrete example, consider the linear case in (\ref{rmk-1-xi-phi}) with the matrix 
  \begin{align}
  a \equiv \begin{pmatrix}
    I_{m\times m} & 0 \\
    0 & \frac{1}{\delta} I_{(n-m)\times (n-m)}
  \end{pmatrix}
  \,,
  \label{rmk-2-mat-a}
\end{align}
where we have recovered the SDE (\ref{linear-case-2}) in the Introduction.
In this case, in analogy to Remark~\ref{rmk-1}, we have $C_1 =0$ in
(\ref{pathwise-bound-case-2}) and the pathwise estimate becomes 
  \begin{align}
    \mathbf{E}\Big(\sup_{0 \le s \le t}\big|\xi(x(s)) - z(s)\big|^2\Big) \le
    \frac{C_2 \delta\,t}{\rho_0^2} e^{Lt}\,.
  \label{pathwise-bound-case-2-linear}
\end{align}
In the general case, however, it is important to point out that the error bound (\ref{pathwise-bound-case-2}) 
 can still be large even when the time scale separation parameter $\delta$  is small.
 We refer to Remark~\ref{rmk-4}, as well as the previous work~\cite{liu2010} for relevant discussions when $\xi$ is linear and $a$ is non-identity matrix. 
   \label{rmk-2}
\end{remark}
\vspace{0.1cm}

  \textbf{Case $3$}. In the third case, we consider the combination of the
  above two cases, i.e., the potential $V$ is given in (\ref{v-v0-v1}) and the
  matrix $a$ has large eigenvalues such that (\ref{pi-a-order-delta}) is satisfied. 
  The following lemma is a direct application of Lemma~\ref{lemma-spectral-gap-stiff-v} and
  Lemma~\ref{lemma-spectral-gap-case-2}. 
  \begin{lemma}
    Given $\epsilon, \delta > 0$. Under the assumptions of
    Lemma~\ref{lemma-spectral-gap-stiff-v} and
assume the matrices $a$ and $\Pi$ satisfy (\ref{a-pi-elliptic-eps})
	for some $c_2 > 0$, which is independent of $\epsilon, \delta$. Then 
	$\exists \epsilon_0,\,C \ge 0$, which may depend on $V_0$, $a$, $\xi$ and $\beta$, such that when $\epsilon \le \epsilon_0$, we have the Poincar{\'e} inequality 
       \begin{align}
	 \mbox{Var}_{\mu_z}(f) = \int_{\Sigma_z} f^2\,d\mu_z - 
	 \Big(\int_{\Sigma_z} f\,d\mu_z\Big)^2 \le \frac{C \epsilon\delta}{c_2 K} \mathcal{E}_z(f,f)\,,
	 \label{poincare-ineq-case-3}
       \end{align}
        for all functions $f : \Sigma_z\rightarrow \mathbb{R}$ which satisfy $\mathcal{E}_z(f,f) < +\infty$.
	\label{lemma-spectral-gap-case-3}
  \end{lemma}
  \begin{proof}
    In the proof of Lemma~\ref{lemma-spectral-gap-stiff-v} in Appendix~\ref{sec-app-2},
    we have actually proved that (\ref{poincare-ineq-no-eps-case-2}) is satisfied 
    with $\rho_0 = \frac{K}{C \epsilon}$, for some constant $C>0$. See (\ref{nu-to-mu}) for details. 
    Therefore, the Poincar{\'e} inequality (\ref{poincare-ineq-case-3})
    follows as a direct consequence of Lemma~\ref{lemma-spectral-gap-case-2}.
  \end{proof}

Applying Theorem~\ref{thm-pathwise-estimate} and
Lemma~\ref{lemma-spectral-gap-case-3}, we can obtain the pathwise estimate of
the effective dynamics in this case.
\vspace{0.1cm}

\begin{corollary}
Assume Assumptions~\ref{assump-0}-\ref{assump-1} hold. Let the potential $V$ be given in
(\ref{v-v0-v1}), where $\epsilon>0$ is a small parameter, and suppose that the
condition (\ref{decomp-pi-orthgonal-v1}) is satisfied.
 Assume the assumptions of Lemma~\ref{lemma-spectral-gap-stiff-v} hold
  uniformly for $z \in \mathbb{R}^m$.
Matrices $a$ and $\Pi$ satisfy (\ref{a-pi-elliptic-eps-upper-lower})
	for some $0 < c_2 \le c_2'$, which are independent of $\epsilon, \delta$.
  $x(s)$ satisfies the SDE (\ref{dynamics-a}) starting from
  $x(0)\sim \mu$, and $z(s)$ is the effective dynamics (\ref{modified-eff-dynamics}) with $z(0) =
  \xi(x(0))$.
  Define $L = 3L_b^2+\frac{48L_\sigma^2}{\beta}+1$.
Then $\exists \epsilon_0 \ge 0$, such that when $\epsilon \le \epsilon_0$, we have 
\begin{align}
  \mathbf{E}\Big(\sup_{0 \le s \le t}\big|\xi(x(s)) - z(s)\big|^2\Big) \le 
    t\Big(\frac{C_1\epsilon}{K} + \frac{C_2\epsilon^2 \delta}{K^2}\Big) e^{Lt}\,,
  \label{pathwise-bound-case-3}
\end{align}
for some constants $C_1, C_2>0$ which are independent of $\epsilon$, $\delta$ and $K$.
    \label{corollary-case-3}
\end{corollary}
\begin{proof}
  The proof is similar to that of Corollary~\ref{corollary-case-2}, 
       by noticing that $\kappa_1^2, \kappa_2^2 \le \frac{C}{\delta}$ and
      Assumption~\ref{assump-4} holds with $\rho=\frac{C K}{\epsilon\delta}$, for some $C>0$. 
\end{proof}
  \begin{remark}
    Consider the linear case (\ref{rmk-1-xi-phi}) in Remark~\ref{rmk-1}
    and assume the potential $V$ is given in (\ref{v-v0-v1}) with $V_1(x)=V_1(x_{m+1}, \cdots, x_n)$. Also let the matrix $a$ be given in (\ref{rmk-2-mat-a}), then we recover the SDE (\ref{linear-case-3}) in the Introduction.
    Correspondingly, in this case $C_1=0$ in (\ref{pathwise-bound-case-3}) and
    therefore we have the pathwise estimate
 \begin{align}
  \mathbf{E}\Big(\sup_{0 \le s \le t}\big|\xi(x(s)) - z(s)\big|^2\Big) \le 
    \frac{C_2\epsilon^2 \delta\,t}{K^2} e^{Lt}\,.
  \label{pathwise-bound-case-3-linear}
\end{align}
       \label{rmk-3}
  \end{remark}
\vspace{0.1cm}

\section{Preliminary pathwise estimates : Proof of
Proposition~\ref{prop-pathwise-1}}
\label{sec-preliminary-estimate}
In this section, after deriving two useful lemmas, we prove our first version of pathwise estimate of the effective dynamics. 
\begin{lemma}
  Recall that $\mu$ is the invariant measure in (\ref{mu-density}) and let
  $\varphi$ be the function defined in (\ref{f-lxi-b-diff}). Under
Assumption~\ref{assump-0}, Assumption~\ref{assump-2} and Assumption~\ref{assump-4}, we have 
  \begin{align}
    \mathbf{E}_\mu |\varphi|^2 = \int_{\mathbb{R}^n} |\varphi|^2\,d\mu \le \frac{\kappa_1^2}{\beta\rho}\,.
  \end{align}
  \label{lemma-f2-by-poincare}
\end{lemma}
\begin{proof}
  From the definition of the function $\widetilde{b}$ in (\ref{b-sigma-eff}), we have
  $\int_{\Sigma_z} \varphi(x)\,d\mu_z(x) = 0$, $\forall z \in \mathbb{R}^m$.
  Furthermore, (\ref{mat-pi}) and (\ref{f-lxi-b-diff}) imply that $\Pi\nabla
  \varphi_i = \Pi\nabla \mathcal{L}\xi_i$, for $1 \le i \le m$. 
  Therefore, applying Assumptions~\ref{assump-2}-\ref{assump-4} and
  using the expression of the Dirichlet form $\mathcal{E}_z$ in (\ref{dirichlet-form-z}), we can derive 
  \begin{align*}
    \int_{\mathbb{R}^n} |\varphi|^2\, d\mu 
    =& \sum_{i=1}^m \int_{\mathbb{R}^m} \Big(\int_{\Sigma_z} |\varphi_i|^2\,d\mu_z\Big)\,Q(z)\,dz\\
    \le & \sum_{i=1}^m\frac{1}{\rho} \int_{\mathbb{R}^m}
    \mathcal{E}_z(\varphi_i, \varphi_i)\, Q(z)\,dz\\ 
  =& \sum_{i=1}^m\frac{1}{\beta \rho} \int_{\mathbb{R}^n} (\Pi \nabla \varphi_i)\cdot (a\Pi
    \nabla \varphi_i)\, d\mu \\
    =& \sum_{i=1}^m\frac{1}{\beta \rho} \int_{\mathbb{R}^n} (\Pi \nabla
    \mathcal{L}\xi_i)\cdot (a\Pi \nabla \mathcal{L}\xi_i)\, d\mu =
    \frac{\kappa_1^2}{\beta \rho}\,.
\end{align*}
\end{proof}

We also need to estimate the Frobenius norm of the difference of the two matrices which appeared in the
noise term of the equation (\ref{xi-z-diff}).
\begin{lemma}
  Assume Assumption~\ref{assump-0} holds.
  Recall the positive definite symmetric matrix functions $\widetilde{\sigma}$, $A$
  defined in (\ref{b-sigma-eff}) and (\ref{mat-a}), respectively.
  We have 
  \begin{align}
\mathbf{E}_\mu \big\|A-\widetilde{\sigma}\circ\xi\big\|^2_F =  \mathbf{E}_\mu \big\|A - (\mathbf{E}_{\mu_z}A)\circ\xi\big\|_F^2 +
      \mathbf{E}_\mu \big\| \big(\widetilde{\sigma} -
      \mathbf{E}_{\mu_z}A\big)\circ \xi\big\|_F^2 \,,
      \label{frobenius-equality}
  \end{align}
and 
  \begin{align}
    \mathbf{E}_\mu \big\|A- (\mathbf{E}_{\mu_z}A)\circ\xi\big\|^2_F \le 
    \mathbf{E}_\mu \big\|A-\widetilde{\sigma}\circ\xi\big\|^2_F \le 2\,\mathbf{E}_\mu \big\|A- (\mathbf{E}_{\mu_z}A)\circ\xi\big\|^2_F \,.
      \label{frobenius-inequality-1}
  \end{align}
  Further suppose that Assumption~\ref{assump-2} and Assumption~\ref{assump-4}
  hold, then we have 
  \begin{align}
    \mathbf{E}_\mu \big\|A-\widetilde{\sigma}\circ\xi\big\|^2_F \le
    \frac{2\kappa_2^2}{\beta \rho}\,.
      \label{frobenius-inequality-2}
  \end{align}
  \label{lemma-Frobenius-estimate}
\end{lemma}
\begin{proof}
  From the definitions (\ref{b-sigma-eff}) and (\ref{mat-a}), we
  have $\widetilde{\sigma}(z) = (\mathbf{E}_{\mu_z}A^2)^{\frac{1}{2}}$,
  $\forall z \in \mathbb{R}^m$. Direct calculation shows that 
  \begin{align}
    \begin{split}
    &\mathbf{E}_\mu \big\|A-\widetilde{\sigma}\circ\xi\big\|^2_F\\
    =& \int_{\mathbb{R}^m} \Big( \mathbf{E}_{\mu_z} \big\|A-
    (\mathbf{E}_{\mu_z}A^2)^{\frac{1}{2}}\big\|^2_F\Big) Q(z)\, dz\\
    =& \int_{\mathbb{R}^m} \Big[ 2\,\mbox{tr}(\mathbf{E}_{\mu_z}A^2) -
      2\,\mbox{tr}\Big((\mathbf{E}_{\mu_z}A)\,(\mathbf{E}_{\mu_z}A^2)^{\frac{1}{2}}\Big)
      \Big] Q(z)\, dz\\
      =&\int_{\mathbb{R}^m} \Big[ \,\mathbf{E}_{\mu_z}\Big(\mbox{tr}\big(A - \mathbf{E}_{\mu_z}A)^2\Big) + 
      \mbox{tr}\Big(\big((\mathbf{E}_{\mu_z}A^2)^{\frac{1}{2}} -
      \mathbf{E}_{\mu_z}A\big)^2\Big) \Big] Q(z)\, dz\\
      = & \mathbf{E}_\mu \big\|A - (\mathbf{E}_{\mu_z}A)\circ\xi\big\|_F^2 +
      \mathbf{E}_\mu \big\| \big(\widetilde{\sigma} -
      \mathbf{E}_{\mu_z}A\big)\circ \xi\big\|_F^2 \,,
    \end{split}
    \label{lemma-Frobenius-proof-ineq1}
  \end{align}
  from which the equality (\ref{frobenius-equality}) and the lower bound in (\ref{frobenius-inequality-1}) follow.
  For the upper bound in (\ref{frobenius-inequality-1}), applying Lieb's concavity
  theorem~\cite{LIEB1973,ANDO1979}, we can estimate
  \begin{align}
    \mbox{tr}\Big((\mathbf{E}_{\mu_z}A)(\mathbf{E}_{\mu_z}A^2)^{\frac{1}{2}}\Big)
    \ge \mathbf{E}_{\mu_z} \mbox{tr}\Big((\mathbf{E}_{\mu_z}A)A(\cdot)\Big) = 
    \mbox{tr}\Big(\big(\mathbf{E}_{\mu_z} A\big)^2\Big)\,,
  \end{align}
  and therefore (\ref{lemma-Frobenius-proof-ineq1}) implies
  \begin{align}
    \begin{split}
    &\mathbf{E}_\mu \big\|A-\widetilde{\sigma}\circ\xi\big\|^2_F\\
      \le &\, 2\int_{\mathbb{R}^m} \Big[
      \,\mbox{tr}\big(\mathbf{E}_{\mu_z}A^2\big) - \,
      \mbox{tr}\Big(\big(\mathbf{E}_{\mu_z}A)^2\Big) \Big] Q(z)\, dz \\
      = &\, 2\int_{\mathbb{R}^m} \Big[
	\,\mbox{tr}\Big(\mathbf{E}_{\mu_z}\big(A-\mathbf{E}_{\mu_z}A\big)^2\Big) \Big] Q(z)\, dz \\
      = &\, 2\,\mathbf{E}_\mu \big\|A- (\mathbf{E}_{\mu_z}A)\circ\xi\big\|^2_F \,.
    \end{split}
    \label{lemma-Frobenius-proof-ineq2}
  \end{align}
  Finally, under Assumption~\ref{assump-2} and Assumption~\ref{assump-4}, applying
Poincar{\'e} inequality, we obtain
\begin{align*}
    &\mathbf{E}_\mu \big\|A-\widetilde{\sigma}\circ\xi\big\|^2_F\\
  \le &\,2\int_{\mathbb{R}^m}\Big[
\mbox{tr}\Big(\mathbf{E}_{\mu_z}\big(A-\mathbf{E}_{\mu_z}A\big)^2\Big) \Big] Q(z)\, dz \\
= &\,2\sum_{1 \le i,j \le m} \int_{\mathbb{R}^m}
\Big[\mathbf{E}_{\mu_z}\big(A_{ij}-\mathbf{E}_{\mu_z}A_{ij}\big)^2\Big]\, Q(z)\, dz \\
\le &\,\frac{2}{\rho}\sum_{1 \le i,j \le m} \int_{\mathbb{R}^m}
\mathcal{E}_z(A_{ij}, A_{ij}) \,Q(z)\, dz =  \frac{2\kappa_2^2}{\beta\rho}\,.
\end{align*}
\end{proof}

Applying the above two lemmas, we are ready to prove the first pathwise estimate
Proposition~\ref{prop-pathwise-1}.
\begin{proof}[Proof of Proposition~\ref{prop-pathwise-1}]
  Recall the function $\varphi$ defined in (\ref{f-lxi-b-diff}). From (\ref{xi-z-diff}), we have 
\begin{align}
  \begin{split}
    \xi(x(t)) - z(t) =& \int_0^t \varphi(x(s))\,ds + \int_0^t
    \Big(\widetilde{b}\big(\xi(x(s))\big) - \widetilde{b}(z(s))\Big)\,ds +
    \sqrt{2\beta^{-1}} M(t)\,,
  \end{split}
  \label{xi-z-diff-sqr}
\end{align}
where $M(t)$ denotes the martingale term
\begin{align}
  M(t) = \int_0^t \Big(A(x(s)) - \widetilde{\sigma}(z(s))\Big)\,d\widetilde{w}(s)\,,
  \label{martingale-m}
\end{align}
with the matrix-valued function $A$ defined in (\ref{mat-a}).
  Therefore, squaring both sides of (\ref{xi-z-diff-sqr}), using
  Assumption~\ref{assump-1} and the elementary inequality
  $(a+b+c)^2 \le 3(a^2 + b^2+c^2)$, $\forall\,a,b,c \in \mathbb{R}$, we obtain
\begin{align*}
  \begin{split}
    \big|\xi(x(t)) - z(t)\big|^2 \le &
3\bigg|\int_0^t \varphi(x(s))\,ds\bigg|^2 +
    3L_b^2\bigg(\int_0^t \big|\xi(x(s)) - z(s)\big|\,ds\bigg)^2 + \frac{6}{\beta}
    \big|M(t)\big|^2 \,.
  \end{split}
\end{align*}
  Taking supremum
  followed by mathematical expectation in the above inequality, we get 
\begin{align}
  \begin{split}
    &\mathbf{E}\Big(\sup_{0 \le t' \le t}\big|\xi(x(t')) - z(t')\big|^2\Big) \\
    \le &
    \,3\mathbf{E}\bigg[\sup_{0 \le t' \le t} \Big|\int_0^{t'} \varphi(x(s))\,ds\Big|^2\bigg] +
    3L_b^2\, \mathbf{E}\Big(\int_0^t \big|\xi(x(s)) -
    z(s)\big|\,ds\Big)^2 + \frac{6}{\beta} \mathbf{E}\sup_{0 \le s \le t} \big|M(s)\big|^2\,.
  \end{split}
  \label{prop-pathwise-eq1}
\end{align}

To estimate the right hand side of (\ref{prop-pathwise-eq1}),
we notice that $x(s)\sim \mu$ for $s \ge 0$, since $\mu$ is the invariant measure and $x(0)\sim \mu$. 
For the first term in (\ref{prop-pathwise-eq1}), using Cauchy-Schwarz inequality and Lemma~\ref{lemma-f2-by-poincare}, we have
\begin{align}
  \begin{split}
  &\mathbf{E}\bigg[\sup_{0 \le t' \le t} \Big|\int_0^{t'}
\varphi(x(s))\,ds\Big|^2\bigg] \\
\le&  
\mathbf{E}\bigg[\sup_{0 \le t' \le t} \Big(t' \int_0^{t'}
\big|\varphi(x(s))\big|^2\,ds\Big)\bigg] \\
\le & t \,
\mathbf{E}\bigg[\int_0^{t} \big|\varphi(x(s))\big|^2\,ds\bigg] \\
=& t^2 \int_{\mathbb{R}^n} |\varphi|^2 d\mu \le \frac{\kappa_1^2\,t^2}{\beta \rho} \,.
\end{split}
\label{prop-pathwise-term-1}
\end{align}
The last term in (\ref{prop-pathwise-eq1}) can be estimated by applying Doob's martingale inequality as
\begin{align}
   \mathbf{E}\sup_{0 \le s \le t} \big|M(s)\big|^2 \notag \le&\, 4\mathbf{E} |M(t)|^2 \notag\\
   =&\, 4 \mathbf{E} \int_0^t \big\|A(x(s)) - \widetilde{\sigma}(z(s))\big\|_F^2\,ds\notag\\
  \le\, & 8 \int_0^t \mathbf{E}\big\|A(x(s)) -
  \widetilde{\sigma}\big(\xi(x(s))\big)\big\|_F^2\,ds + 
  8 \int_0^t \mathbf{E}\big\|\widetilde{\sigma}\big(\xi(x(s))\big) - \widetilde{\sigma}(z(s))\big\|_F^2\,ds\notag\\
  \le&\, \frac{16\kappa_2^2 t}{\beta \rho} + 8L_\sigma^2 \int_0^t
  \mathbf{E}\big|\xi(x(s)) -
  z(s)\big|^2\,ds\,, 
  \label{prop-pathwise-term-3}
\end{align}
where we have used Assumption~\ref{assump-1}, Lemma~\ref{lemma-Frobenius-estimate}, together with the fact that $x(s)\sim \mu$. 
Combining (\ref{prop-pathwise-eq1}), (\ref{prop-pathwise-term-1}), and (\ref{prop-pathwise-term-3}), we get
\begin{align}
  \begin{split}
    &\mathbf{E}\Big(\sup_{0 \le t' \le t}\big|\xi(x(t')) - z(t')\big|^2\Big) \\
    \le &
  \Big(\frac{3\kappa_1^2\,t^2}{\beta \rho} + \frac{96\kappa_2^2 t}{\beta^2 \rho}\Big) 
    +3L_b^2\, \mathbf{E}\Big(\int_0^t \sup_{0 \le t' \le s}
    \big|\xi(x(t')) - z(t')\big|\,ds\Big)^2
  + \frac{48L_\sigma^2}{\beta} \int_0^t \Big(\mathbf{E}\sup_{0 \le t' \le s}
  \big|\xi(x(t')) - z(t')\big|^2\Big)\,ds \,.
  \end{split}
  \label{prop-pathwise-pre-gronwall}
\end{align}
  The conclusion follows by applying Lemma~\ref{lemma-gronwall-like-ineq} below.
\end{proof}
\begin{remark}
  Applying Cauchy-Schwarz inequality to the right hand side of
  (\ref{prop-pathwise-pre-gronwall}), we can get
\begin{align}
  \begin{split}
    &\mathbf{E}\Big(\sup_{0 \le t' \le t}\big|\xi(x(t')) - z(t')\big|^2\Big) \\
    \le &
  \Big(\frac{3\kappa_1^2\,t^2}{\beta \rho} + \frac{96\kappa_2^2 t}{\beta^2 \rho}\Big) 
   +\Big(3L_b^2\,t+ \frac{48L_\sigma^2}{\beta}\Big)
    \, \int_0^t \Big(\mathbf{E}\sup_{0 \le t' \le s} \big|\xi(x(t')) -
    z(t')\big|^2\Big)\,ds\,,
  \end{split}
\end{align}
and therefore Gronwall's inequality directly implies 
  \begin{align*}
\mathbf{E}\Big(\sup_{0 \le t' \le t}\big|\xi(x(t')) - z(t')\big|^2\Big)
\le 
  \Big(\frac{3\kappa_1^2\,t^2}{\beta \rho} + \frac{96\kappa_2^2 t}{\beta^2
    \rho}\Big) e^{Lt}\,,
  \end{align*}
  with $L = 3L_b^2\,t+ \frac{48L_\sigma^2}{\beta}$. 
  Notice that, however, using this argument the constant $L$ will depend on the time $t$.
  Instead, Lemma~\ref{lemma-gronwall-like-ineq} below allows us to obtain
  an upper bound where the constant $L= 3L_b^2+
  \frac{48L_\sigma^2}{\beta}+1$, which is independent of $t$.
\end{remark}
\begin{lemma}
  Let $f(t) \in \mathbb{R}$ be a function on $t\in [0, +\infty)$ taking random values, such that $\mathbf{E}
  (f(t))^2 < +\infty$ for all $t\ge 0$.
  Further assume that $f$ satisfies the inequality 
  \begin{align}
    \mathbf{E}\big(f(t)\big)^2 \le g(t) + C_1 \mathbf{E}\Big(\int_0^t f(s) ds\Big)^2 + C_2
    \int_0^t \mathbf{E} \big(f(s)\big)^2 ds\,, \quad \forall\, t \ge 0\,,
    \label{eqn-1-gronwall-like-lemma}
  \end{align}
  where the constants $C_1, C_2\ge 0$, and $g\ge 0$ is a function of $t\in [0, +\infty)$. We have 
  \begin{align}
    \mathbf{E}\big(f(t)\big)^2 \le g(t) + (C_1 + C_2) \int_0^t e^{(C_1 +C_2 +
    1)(t-s)} g(s)\,ds\,,\quad \forall t \ge 0\,.
  \label{eqn-fsquare-lemma-gronwall-like-ineq}
  \end{align}
  In particular, when the function $g$ is non-decreasing, we have
  \begin{align}
    \mathbf{E}\big(f(t)\big)^2 \le g(t)\,e^{(C_1 + C_2 + 1) t}\,,\quad \forall t \ge 0\,.
  \label{eqn-fsquare-g-increase-lemma-gronwall-like-ineq}
  \end{align}
  \label{lemma-gronwall-like-ineq}
\end{lemma}
\begin{proof}
  For $t \ge 0$, we define the function 
  \begin{align}
    F(t)= C_1 \mathbf{E}\Big(\int_0^t f(s) ds\Big)^2 + C_2 \int_0^t \mathbf{E} \big(f(s)\big)^2 ds\,.
    \label{eqn-2-gronwall-like-lemma}
  \end{align}
From (\ref{eqn-1-gronwall-like-lemma}) and (\ref{eqn-2-gronwall-like-lemma}), we can compute 
\begin{align*}
  \frac{dF}{dt} =&\, 2 C_1 \mathbf{E}\Big(f(t) \int_0^t f(s) ds\Big) + C_2
  \mathbf{E} \big(f(t)\big)^2\\
  \le&\, C_1 \mathbf{E}\Big(\int_0^t f(s) ds\Big)^2 + (C_1 + C_2)
  \mathbf{E} \big(f(t)\big)^2 \\
  \le&\, F(t) + (C_1 + C_2) \big(g(t) + F(t)\big)\\
  = &\, (C_1 + C_2 + 1) F(t) + (C_1 + C_2) g(t)\,. 
\end{align*}
Since $F(0) = 0$, after integration we obtain 
\begin{align}
  F(t) \le (C_1 + C_2) \int_0^t e^{(C_1 +C_2 + 1)(t-s)} g(s)\,ds\,.
  \label{bound-f-gronwall-lemma}
\end{align}
  The inequalities (\ref{eqn-fsquare-lemma-gronwall-like-ineq}) and (\ref{eqn-fsquare-g-increase-lemma-gronwall-like-ineq})
  follow by substituting (\ref{bound-f-gronwall-lemma}) into (\ref{eqn-1-gronwall-like-lemma}). 
\end{proof}
\section{Pathwise estimates of effective dynamics: Proof of Theorem~\ref{thm-pathwise-estimate} and Theorem~\ref{thm-pathwise-estimate-mu-bar}} 
  \label{sec-pathwise-approx}
In this section, we prove the Theorem~\ref{thm-pathwise-estimate} and Theorem~\ref{thm-pathwise-estimate-mu-bar},
which improve the pathwise estimate result in Proposition~\ref{prop-pathwise-1}.
The main tool we will use is the forward-backward martingale approach developed in~\cite{lyons1994,komorowski2012fluctuations}.

Following the argument in~\cite{LEGOLL2017-pathwise}, we first establish a technique
result. Given $\Psi \in C^2(\mathbb{R}^n, \mathbb{R}^n)$, we denote by
  \begin{align}
 \nabla^*\Psi=\beta \nabla V \cdot \Psi - \mbox{div}\Psi
 \label{nabla-adjoint}
 \end{align}
  the adjoint of the gradient operator $\nabla$ with
  respect to the probability measure $\mu$, i.e.,
  \begin{align}
    \int_{\mathbb{R}^n} h\nabla^*\Psi\,d\mu = \int_{\mathbb{R}^n} \Psi\cdot \nabla h\, d\mu\,
  \end{align}
  holds for any $C^1$ function $h :\mathbb{R}^n \rightarrow \mathbb{R}$. 
  We have the following estimate.
  \begin{prop}
    Recall that $\mu$ is the invariant measure in (\ref{mu-density}) and the
    matrix $a$ satisfies the uniform elliptic condition (\ref{a-elliptic}). 
    Let $\Psi \in C^2(\mathbb{R}^n, \mathbb{R}^n)$ and 
    $x(s)$ be the dynamics in (\ref{dynamics-a}) starting from $x(0)\sim \mu$.
    For all $t \ge 0$, we have 
  \begin{align}
    \mathbf{E}\bigg[\sup_{0 \le t' \le t} \Big|\int_0^{t'} \nabla^*\Psi(x(s))\,
    ds\Big|^2\bigg]
    \le \frac{27\beta t}{2} \int_{\mathbb{R}^n} \Psi\cdot (a^{-1} \Psi)\, d\mu\,.
  \end{align}
  \label{prop-forward-backward}
  \end{prop}
  \begin{proof}
    We will only sketch the proof since the argument resembles the one in
    \cite{LEGOLL2017-pathwise}, with minor modifications due to the appearance of the matrix $a$. 
    First of all, condition (\ref{a-elliptic}) implies that the matrix $a$ is
    positive definite and therefore invertible.
    Given $\eta > 0$, we consider the function 
    $\omega_\eta : \mathbb{R}^n \rightarrow \mathbb{R}$ which solves the PDE 
  \begin{align}
    \eta \omega_{\eta} - \mathcal{L} \omega_{\eta} = - \nabla^*\Psi\,.
    \label{eta-w}
  \end{align}
    Multiplying both sides of (\ref{eta-w}) by $\omega_{\eta}$ and
  integrating with respect to $\mu$, we obtain
  \begin{align*}
    \eta \int_{\mathbb{R}^n}  \omega_{\eta}^2\, d\mu +
    \mathcal{E}(\omega_{\eta}, \omega_{\eta}) = - \int_{\mathbb{R}^n} \Psi\cdot \nabla
    \omega_{\eta}\, d\mu \,,
  \end{align*}
    where $\mathcal{E}$ is the Dirichlet form in (\ref{lfg-flg}).
    Applying Cauchy-Schwarz inequality to the right hand side of the above
    equality, we can estimate
    \begin{align*}
      & \eta \int_{\mathbb{R}^n}  \omega_{\eta}^2\, d\mu +
      \mathcal{E}(\omega_{\eta}, \omega_{\eta})  \\
      \le&  \Big(\int_{\mathbb{R}^n} \Psi\cdot (a^{-1} \Psi)\, d\mu\Big)^{\frac{1}{2}}
    \Big(\int_{\mathbb{R}^n} \nabla \omega_{\eta} \cdot (a\nabla \omega_{\eta})\,
    d\mu\Big)^{\frac{1}{2}} \\
       = & \Big(\int_{\mathbb{R}^n} \Psi\cdot (a^{-1} \Psi)\, d\mu\Big)^{\frac{1}{2}}
       \Big(\beta \mathcal{E}(\omega_{\eta},\omega_{\eta})\Big)^{\frac{1}{2}}\,,
  \end{align*}
  from which we can deduce that 
  \begin{align}
    \begin{split}
    & \mathcal{E}(\omega_{\eta},\omega_{\eta})\le 
      \beta \int_{\mathbb{R}^n} \Psi\cdot (a^{-1} \Psi)\, d\mu\,, \\
      &  \lim_{\eta \rightarrow 0} \bigg(\eta^2 \int_{\mathbb{R}^n}
      \omega_{\eta}^2\, d\mu\bigg) = 0\,.
  \end{split}
  \label{w-eta-bound}
  \end{align}
   Now let us define the process $y(t') := x(t-t')$ for $t' \in [0, t]$.
    Since $x(s)$ is both reversible and in stationary, the process
    $(y(s))_{\{0 \le s \le t\}}$
    satisfies the same law as $(x(s))_{\{0 \le s \le t\}}$. And we can assume that 
    there is another $n'$-dimensional Brownian motion $\bar{w}$, such that 
    \begin{align}
    dy(s) = -a(y(s))\nabla V(y(s))\, ds + \frac{1}{\beta} (\nabla \cdot
	a)(y(s))\, ds + \sqrt{2\beta^{-1}} \sigma(y(s))\, d\bar{w}(s)\,, \quad
	s \in [0, t]\,.
    \label{dynamics-y}
  \end{align}

  Applying Ito's formula, we have 
  \begin{align}
    \omega_\eta(x(t')) - \omega_\eta(x(0)) = \int_0^{t'} (\mathcal{L}\omega_\eta)(x(s))\, ds +
    \sqrt{2\beta^{-1}} \int_0^{t'} \big[(\nabla
    \omega_\eta)^T\sigma\big](x(s))\,dw(s)\,,
    \label{w-ito-x}
  \end{align}
  and similarly, 
  \begin{align}
    \omega_\eta(y(t)) - \omega_\eta(y(t-t')) = \int_{t-t'}^t (\mathcal{L}\omega_\eta)(y(s))\, ds +
    \sqrt{2\beta^{-1}} \int_{t-t'}^t \big[(\nabla \omega_\eta)^T\sigma\big](y(s))\,
    d\bar{w}(s)\,.
    \label{w-ito-y}
  \end{align}
    Let us denote 
  \begin{align*}
    M(t') = \int_0^{t'} \big[(\nabla \omega_\eta)^T\sigma\big](x(s))\, dw(s)\,,\quad
    \mbox{and}\,\quad \overline{M}(t') =
    \int_0^{t'} \big[(\nabla \omega_\eta)^T\sigma\big](y(s))\, d\bar{w}(s)\,.
  \end{align*}
  Adding up (\ref{w-ito-x}) and (\ref{w-ito-y}), we obtain 
  \begin{align}
    2 \int_0^{t'} (\mathcal{L}\omega_\eta)(x(s))\,ds + \sqrt{2\beta^{-1}} \Big[M(t') +
    \overline{M}(t) - \overline{M}(t-t')\Big]  = 0\,,
    \label{int-lw-m}
  \end{align}
  where the terms $M(t')$, $\overline{M}(t)$ and $\overline{M}(t-t')$ can be
    bounded using Doob's martingale inequality. We refer to~\cite{LEGOLL2017-pathwise} for details.
  Combining these upper bounds with (\ref{int-lw-m}) and (\ref{w-eta-bound}),
  we can obtain
  \begin{align}
    \mathbf{E}\bigg[\sup_{0 \le t' \le t} \Big|\int_0^{t'} (\mathcal{L}\omega_{\eta})(x(s))\,
    ds\Big|^2\bigg]
    \le \frac{27t}{2} \mathcal{E}(\omega_\eta, \omega_\eta)\,\le
    \frac{27\beta t}{2} \int_{\mathbb{R}^n} \Psi\cdot (a^{-1} \Psi)\, d\mu\,.
    \label{path-int-bound}
  \end{align}
  Letting $\eta \rightarrow 0$, using (\ref{eta-w}), (\ref{w-eta-bound})
  and the estimate (\ref{path-int-bound}), we conclude that 
  \begin{align*}
    \mathbf{E}\bigg[\sup_{0 \le t' \le t} \Big|\int_0^{t'} (\nabla^*\Psi)(x(s))\,
    ds\Big|^2\bigg]
    \le \frac{27\beta t}{2} \int_{\mathbb{R}^n} \Psi\cdot (a^{-1} \Psi)\, d\mu\,.
  \end{align*}
\end{proof}

  To proceed, we consider the equation
  \begin{align}
    \mathcal{L}_0 u = \varphi\,,
    \label{poisson-sigma-z}
  \end{align}
  where $u : \mathbb{R}^n \rightarrow \mathbb{R}^m$, $\mathcal{L}_0$ is the
  operator in (\ref{l-0-1}) and  $\varphi$ is the function introduced in (\ref{f-lxi-b-diff}).  For each $z \in
  \mathbb{R}^m$, (\ref{poisson-sigma-z}) can be viewed as a Poisson equation
on the submanifold $\Sigma_z$ for the components of $u$. 
Applying Proposition~\ref{prop-forward-backward}, we can obtain the following result. 
  \begin{lemma}
    Assume Assumption~\ref{assump-0}, Assumption~\ref{assump-2} and Assumption~\ref{assump-4} hold. Let $x(s)$
    be the dynamics in (\ref{dynamics-a}) starting from $x(0)\sim \mu$.
    For all $t \ge 0$, we have 
  \begin{align*}
    \mathbf{E}\bigg[\sup_{0 \le t' \le t} \Big|\int_0^{t'} \varphi(x(s))\, ds\Big|^2\bigg]
     \le \frac{27\kappa_1^2 t}{2\beta\rho^2}\,.
  \end{align*}
    \label{lemma-int-f-sqr2}
  \end{lemma} 
  \begin{proof}
    Using the equation (\ref{poisson-sigma-z}) and the fact that $\mathcal{E}_z$ satisfies
    Poincar{\'e} inequality (Assumption~\ref{assump-4}), we can deduce that 
  \begin{align}
    \mathcal{E}_z(u_i, u_i) \le \frac{1}{\rho} \int_{\Sigma_z} \varphi_i^2\, d\mu_z\,,
    \label{grad-u-bound-by-f-l2}
  \end{align}
  for $1 \le i \le m$ and for all $z \in \mathbb{R}^m$. We refer to~\cite[Lemma
    9]{LEGOLL2017-pathwise} for details. Recalling the matrix $\Pi$ in (\ref{proj-pi}) and (\ref{mat-pi}),
    let us define $\Psi^{(i)} = a \Pi \nabla u_i \in \mathbb{R}^n$, for $1 \le i \le m$. 
    From (\ref{nabla-adjoint}) and the expression of $\mathcal{L}_0$ in
    (\ref{l-0-1}), direct calculation shows that
  \begin{align}
    \nabla^* \Psi^{(i)} 
    = \beta\nabla V\cdot \Psi^{(i)} - \mbox{div}\,\Psi^{(i)} 
    = \beta\nabla V\cdot (a\Pi \nabla u_i) - \mbox{div}(a\Pi \nabla u_i) 
    = -\beta \mathcal{L}_0 u_i \,.
  \end{align}
    Therefore, applying Proposition~\ref{prop-forward-backward}, Lemma~\ref{lemma-f2-by-poincare} and the expression in
  (\ref{dirichlet-form-z}), we can derive
  \begin{align*}
    &\mathbf{E}\bigg[\sup_{0 \le t' \le t} \Big|\int_0^{t'} \varphi(x(s))\, ds\Big|^2\bigg]\\
    =&\,\mathbf{E}\bigg[\sup_{0 \le t' \le t} \sum_{i=1}^m\Big|\int_0^{t'} \varphi_i(x(s))\, ds\Big|^2\bigg]\\
    \le&\,\sum_{i=1}^m\mathbf{E}\bigg[\sup_{0 \le t' \le t} \Big|\int_0^{t'} \varphi_i(x(s))\, ds\Big|^2\bigg]\\
    =&    \sum_{i=1}^m\mathbf{E}\bigg[\sup_{0 \le t' \le t} \Big|\int_0^{t'} (\mathcal{L}_0
    u_i)(x(s))\, ds\Big|^2\bigg] \\
    \le& \frac{27t}{2\beta} \sum_{i=1}^m \int_{\mathbb{R}^n} (a\Pi\nabla
    u_i)\cdot \Big[a^{-1} (a \Pi\nabla u_i)\Big]\, d\mu\\
    =& \frac{27t}{2\beta} \sum_{i=1}^m \int_{\mathbb{R}^m}\bigg[ \int_{\Sigma_z} (a\Pi\nabla u_i)\cdot (\Pi\nabla u_i)\, d\mu_z\bigg] Q(z)\, dz \\
    =& \frac{27t}{2} \sum_{i=1}^m \int_{\mathbb{R}^m} \mathcal{E}_z(u_i, u_i)\, Q(z)\, dz \\
    \le & \frac{27t}{2\rho} \sum_{i=1}^m \int_{\mathbb{R}^m} \Big(\int_{\Sigma_z} \varphi_i^2\, d\mu_z\Big) Q(z)\,
    dz \le \frac{27\kappa_1^2 t}{2\beta\rho^2}\,,
  \end{align*}
  where the inequality (\ref{grad-u-bound-by-f-l2}) has been used.
\end{proof}

We are ready to prove Theorem~\ref{thm-pathwise-estimate}.
\begin{proof}[Proof of Theorem~\ref{thm-pathwise-estimate}]
  The proof is similar to that of Proposition~\ref{prop-pathwise-1} in
  Section~\ref{sec-preliminary-estimate}, with a few modifications. 
  Specifically, in analogy to the inequalities (\ref{prop-pathwise-eq1}) and (\ref{prop-pathwise-term-3}), 
  we have 
\begin{align}
  \begin{split}
    &\mathbf{E}\Big(\sup_{0 \le t' \le t}\big|\xi(x(t')) - z(t')\big|^2\Big) \\
    \le &
    \,3\mathbf{E}\bigg[\sup_{0 \le t' \le t} \Big|\int_0^{t'} \varphi(x(s))\,ds\Big|^2\bigg] +
    3L_b^2\, \mathbf{E}\Big(\int_0^t \big|\xi(x(s)) - z(s)\big|\,ds\Big)^2 + \frac{6}{\beta} \mathbf{E}\sup_{0 \le s \le t} \big|M(s)\big|^2\,.
  \end{split}
  \label{sup-xi-zt-diff-sqr2}
\end{align}
where the last term can be estimated using Doob's martingale inequality as
\begin{align}
   \mathbf{E}\sup_{0 \le s \le t} \big|M(s)\big|^2 
  \le \frac{16\kappa_2^2 t}{\beta \rho} + 8L_\sigma^2 \int_0^t
  \mathbf{E}\big|\xi(x(s)) -
  z(s)\big|^2\,ds\,.\label{mt-doob-bound} 
\end{align}
  Now, the main different step from Proposition~\ref{prop-pathwise-1} is that we will estimate the first term on the right
hand side of (\ref{sup-xi-zt-diff-sqr2}) by applying
  Lemma~\ref{lemma-int-f-sqr2}.
Combining (\ref{sup-xi-zt-diff-sqr2}), (\ref{mt-doob-bound}) and Lemma~\ref{lemma-int-f-sqr2}, we get
\begin{align*}
  &\mathbf{E}\Big(\sup_{0 \le t' \le t}\big|\xi(x(t')) - z(t')\big|^2\Big)\notag \\
    \le& 
\Big(\frac{81\kappa_1^2}{2\beta \rho^2}
  +\frac{96\kappa_2^2}{\beta^2\rho}\Big)t + 
      3L_b^2\, \mathbf{E}\bigg(\int_0^t \sup_{0 \le t' \le s} \big|\xi(x(t')) -
  z(t')\big|\,ds\bigg)^2 \\
  & + \frac{48L_\sigma^2}{\beta}\int_0^t \mathbf{E}\Big(\sup_{0 \le t' \le s}
  \big|\xi(x(t')) - z(t')\big|^2\Big)\,ds\,.
\end{align*}
  The conclusion follows by applying Lemma~\ref{lemma-gronwall-like-ineq}.
\end{proof}
\vspace{0.2cm}

Finally, we apply Theorem~\ref{thm-pathwise-estimate} to prove Theorem~\ref{thm-pathwise-estimate-mu-bar} for more general initial conditions. 
\vspace{0.1cm}

\begin{proof}[Proof of Theorem~\ref{thm-pathwise-estimate-mu-bar}]
  Let us denote by $\mathbf{E}$ and $\mathbf{E}_{x'}$ the shorthands of
  $\mathbf{E}(\cdot\,|\,x(0)\sim \mu)$, $\mathbf{E}(\cdot\,|\,x(0)=x')$,
  i.e., the expectations with respect to the trajectories $(x(s))_{s\ge 0}$ starting from the
  invariant distribution $\mu$ and the state $x'$, respectively. We will
  frequently use the elementary
  inequality $\sqrt{a+b} \le \sqrt{a} + \sqrt{b}$, $\forall a, b\ge 0$. 
  \begin{enumerate}
    \item
  Applying Theorem~\ref{thm-pathwise-estimate} and Cauchy-Schwarz inequality, we can compute 
  \begin{align*}
    &\mathbf{E}\bigg(\sup_{0 \le s \le t}\big|\xi(x(s)) - z(s)\big|~\Big|~x(0)
  \sim \bar{\mu}\bigg)  \\
  =& \mathbf{E}\bigg[\Big(\sup_{0 \le s \le t}\big|\xi(x(s)) -
z(s)\big|\Big)\,\frac{d\bar{\mu}}{d\mu} \bigg]  \\
\le & 
\bigg[\mathbf{E}\Big(\sup_{0 \le s \le t}\big|\xi(x(s)) -
  z(s)\big|^2\Big)\bigg]^{\frac{1}{2}}\,
  \Big[\int_{\mathbb{R}^n} \Big(\frac{d\bar{\mu}}{d\mu}\Big)^{2} d\mu\Big]^{\frac{1}{2}} \\  
  \le & 
  \sqrt{t}
  \Big(
  \frac{9\kappa_1}{\sqrt{2\beta}\rho} +
  \frac{12\kappa_2}{\beta\sqrt{\rho}}\, \Big) 
  \Big[\int_{\mathbb{R}^n} \Big(\frac{d\bar{\mu}}{d\mu}\Big)^{2} d\mu\Big]^{\frac{1}{2}} 
e^{Lt}\,,
\end{align*}
            where $L =
      \frac{3}{2}L_b^2+\frac{24L_\sigma^2}{\beta}+\frac{1}{2}$.
    \item
       Let us define $h(s) = \int_{\mathbb{R}^n} (p_s-1)^2\,d\mu $ for $s>0$. 
      From the study of the heat kernel estimate~\cite{davies1990heat,weighted-nash}, it is known that $h(s)$ is finite for
      $\forall s>0$. Since $\int_{\mathbb{R}^n} p_s\, d\mu=1$, we have $h(s) =
      \int_{\mathbb{R}^n} p_s^2\,d\mu -1$. Using the Poincar{\'e} inequality
      (\ref{poincare-ineq-xs}) and noticing that $p_s$ satisfies the Kolmogorov equation, we can calculate 
      \begin{align*}
      h'(s) =
      2\int_{\mathbb{R}^n} p_s\,\mathcal{L}p_s\,d\mu = - 2\mathcal{E}(p_s, p_s)
	\le - 2\alpha h(s)\,, 
      \end{align*}
      which implies 
      \begin{align}
      \int_{\mathbb{R}^n} p_{t_1}^2\,d\mu \le 1 + \Big(\int_{\mathbb{R}^n}
	p_{t_0}^2\,d\mu - 1\Big) e^{-2\alpha (t_1 - t_0)} < 1 + e^{-2\alpha (t_1 - t_0)} \int_{\mathbb{R}^n} p_{t_0}^2\,d\mu\,,
	\label{pt-sqr-ineq}
      \end{align}
      for any $0< t_0 \le t_1 \le t$.
      We also introduce the auxiliary process $\bar{z}(s)$, which is the effective dynamics
      (\ref{modified-eff-dynamics}) on $s \in [t_1, t]$, starting from $\bar{z}(t_1) = \xi(x(t_1))$.
      Clearly, we have 
      \begin{align*}
	&\mathbf{E}_{x'}\Big(\sup_{0 \le s \le t}\big|\xi(x(s)) - z(s)\big|\Big)  \\
	\le& \max \bigg\{\mathbf{E}_{x'}\Big(\sup_{0 \le s \le
      t_1}\big|\xi(x(s)) -
	z(s)\big|\Big),\, \mathbf{E}_{x'}\Big(\sup_{t_1 \le s \le
	t}\big|\xi(x(s)) - \bar{z}(s)\big|\Big) + \mathbf{E}_{x'}\Big(\sup_{t_1 \le s \le t}\big|z(s) - \bar{z}(s)\big|\Big)
	\bigg\}\,.
	\end{align*}
	On the time interval $[t_1, t]$, from the estimate
	(\ref{pathwise-bound-mu-bar}) in the previous conclusion and the estimate (\ref{pt-sqr-ineq}), we know 
	\begin{align}
	  \begin{split}
	    & \mathbf{E}_{x'}\Big(\sup_{t_1 \le s \le t}\big|\xi(x(s)) -
	  \bar{z}(s)\big|\Big)\\
	  =& 
	  \mathbf{E}\Big(\sup_{t_1 \le s \le t}\big|\xi(x(s)) - \bar{z}(s)\big|~\Big|~ x(t_1) \sim \mu_{t_1} \Big)\\
	  \le & 
  \sqrt{t-t_1}
  \Big(
  \frac{9\kappa_1}{\sqrt{2\beta}\rho} +
  \frac{12\kappa_2}{\beta\sqrt{\rho}}\, \Big) 
	    \bigg[1 + e^{-\alpha (t_1 - t_0)} \Big(\int_{\mathbb{R}^n} p^2_{t_0}
	    d\mu\Big)^{\frac{1}{2}}\bigg]
	    e^{L(t-t_1)}\,,
	  \end{split}
	  \label{corollary-fixed-initial-ineq-1}
	\end{align}
	where $L =
	\frac{3}{2}L_b^2+\frac{24L_\sigma^2}{\beta}+\frac{1}{2}$.
	Meanwhile, using the same argument as in Proposition~\ref{prop-pathwise-1} and Theorem~\ref{thm-pathwise-estimate}, we can
	obtain the estimate 
	\begin{align}
	  \begin{split}
	     \mathbf{E}_{x'}\Big(\sup_{t_1 \le s \le t}\big|z(s) -
	    \bar{z}(s)\big|^2\Big) 
	     \le &\, 3 \mathbf{E}_{x'}\Big(\big|z(t_1) - \bar{z}(t_1)\big|^2\Big) e^{L_1(t-t_1)} \\
	    \le &\, 3 \mathbf{E}_{x'}\Big(\sup_{0 \le s \le t_1} \big|\xi(x(s)) -
	    z(s)\big|^2\Big) e^{L_1(t-t_1)} \,,
	  \end{split}
	  \label{corollary-fixed-initial-ineq-2}
	\end{align}
	where $L_1 = 3L_b^2 + \frac{24L_\sigma^2}{\beta}+1$, and we have used
	the fact that $\bar{z}(t_1) = \xi(x(t_1))$.

	On the time interval $s \in [0, t_1]$, in analogy to
	(\ref{prop-pathwise-eq1}) in the proof of Proposition~\ref{prop-pathwise-1}, we can obtain 
\begin{align}
  \begin{split}
    &\mathbf{E}_{x'}\Big(\sup_{0 \le t' \le s}\big|\xi(x(t')) - z(t')\big|^2\Big) \\
    \le &
    \,3\mathbf{E}_{x'}\bigg[\sup_{0 \le t' \le s} \Big|\int_0^{t'} \varphi(x(r))\,dr\Big|^2\bigg] +
    3L_b^2\,
    \mathbf{E}_{x'} \bigg(\int_0^s \sup_{0 \le t' \le r} \big|\xi(x(t')) -
    z(t')\big|\,dr\bigg)^2 \\
    & + \frac{6}{\beta} \mathbf{E}_{x'}\Big(\sup_{0 \le t'
    \le s} \big|M(t')\big|^2\Big)\,,  
  \end{split}
\end{align}
where $M(t')$ is the martingale given in (\ref{martingale-m}). Since both
$\varphi$ and $A$ are bounded, applying Doob's martingale inequality, it follows that 
\begin{align*}
  &\mathbf{E}_{x'}\Big(\sup_{0 \le t' \le s}\big|\xi(x(t')) - z(t')\big|^2\Big) \\
    \le &
    \,3C_1^2 s^2 + \frac{96}{\beta} C_2^2\,s
    + 3L_b^2\, \mathbf{E}_{x'} \bigg(\int_0^s \sup_{0 \le t'
    \le r} \big|\xi(x(t')) - z(t')\big|\,dr \bigg)^2 \,,
\end{align*}
which, from Lemma~\ref{lemma-gronwall-like-ineq}, implies 
\begin{align}
  &\mathbf{E}_{x'}\Big(\sup_{0 \le t' \le t_1}\big|\xi(x(t')) - z(t')\big|^2\Big)
    \le 
    \,\Big(3C_1^2 t_1^2 + \frac{96}{\beta} C_2^2\,t_1\Big)
    e^{(3L_b^2 + 1)t_1}\,.
	  \label{corollary-fixed-initial-ineq-3}
  \end{align}
  Combining (\ref{corollary-fixed-initial-ineq-1}),
  (\ref{corollary-fixed-initial-ineq-2}), and (\ref{corollary-fixed-initial-ineq-3}), 
  we conclude that 
  \begin{align*}
    &\mathbf{E}_{x'}\Big(\sup_{0 \le s \le t}\big|\xi(x(s)) - z(s)\big|\Big)  \\
	\le& 
  \sqrt{t}
  \Big(
  \frac{9\kappa_1}{\sqrt{2\beta}\rho} +
  \frac{12\kappa_2}{\beta\sqrt{\rho}}\, \Big) 
	    \bigg[1 + e^{-\alpha (t_1 - t_0)} \Big(\int_{\mathbb{R}^n} p^2_{t_0}
	    d\mu\Big)^{\frac{1}{2}}\bigg] e^{L(t-t_1)} \\
	    & + \sqrt{t_1}\Big(3C_1 \sqrt{t_1} +
	    \frac{18C_2}{\sqrt{\beta}}\Big)
	    e^{(\frac{3}{2}L_b^2+\frac{1}{2}) t_1 +
	    \frac{1}{2} L_1 (t-t_1)} \\
	    \le & 
    \bigg\{\sqrt{t}
  \Big(
  \frac{9\kappa_1}{\sqrt{2\beta}\rho} +
  \frac{12\kappa_2}{\beta\sqrt{\rho}}\, \Big) 
	    \bigg[1 + e^{-\alpha (t_1 - t_0)} \Big(\int_{\mathbb{R}^n} p^2_{t_0}
	    d\mu\Big)^{\frac{1}{2}}\bigg]  
	     + \sqrt{t_1}\Big(3C_1 \sqrt{t_1} +
	     \frac{18C_2}{\sqrt{\beta}}\Big)\bigg\}  e^{Lt} \,.
  \end{align*}
	\end{enumerate}
\end{proof}

\section*{Acknowledgement}
The work of T. Leli{\`e}vre is supported by the European Research Council
under the European Union's Seventh Framework Programme (FP/2007-2013) / ERC
Grant Agreement number 614492.
The work of W. Zhang is supported by the Einstein Foundation Berlin of the
Einstein Center for Mathematics (ECMath) through project CH21. 
Part of the work was done while both authors were attending the program ``Complex High-Dimensional Energy Landscapes'' at IPAM (UCLA), 2017. 
The authors thank the institute for hospitality and support.

  \appendix 
  \section{Coordinate transformation : from nonlinear to linear
  reaction coordinate}
  \label{sec-app-1}
  In this appendix, given a (nonlinear) reaction coordinate function $\xi :
  \mathbb{R}^n \rightarrow \mathbb{R}^m$, we study the coordinate transformation under which the
  original reaction coordinate becomes the mapping onto the first $m$
  components of system's state, i.e., the linear reaction coordinate. 
  Specifically, given $x \in \mathbb{R}^n$, we consider the existence of a function $\phi
  : \Omega_x  \rightarrow \mathbb{R}^{n-m}$, where $\Omega_x \subseteq
  \mathbb{R}^n$ is a neighborhood
  of $x$, such that the map 
 \begin{align}
G(x) = (\xi(x), \phi(x))
   \label{map-g}
 \end{align}
 is one to one from $\Omega_x$ to $\mbox{Im}(G)$.  We further impose that 
 \vspace{0.2cm}
 \begin{align}
 (\nabla\xi a \nabla\phi^T) \equiv 0,\,  \quad \Longleftrightarrow \quad
    (a\nabla\xi_i)_{l} 
    \frac{\partial\phi_{j}}{\partial x_{l}}  = 0\,, \qquad \forall\,~1 \le i\le m\,,\quad 1 \le j \le n-m\,.
    \label{assump-xi-phi}
\end{align} 

Notice that, in this appendix we will adopt the Einstein summation
 convention, i.e., repeated indices indicate summation over a set of indexed
 terms. Recalling the dynamics $x(s)$ in (\ref{dynamics-a}) and applying Ito's
 formula, we know that the dynamics of $\bar{z}(s)=\xi(x(s))$ and
 $\bar{y}(s)=\phi(x(s))$ are given by  
  \begin{align}
    \begin{split}
      d\bar{z}_i(s) =& \big(\mathcal{L}\xi_i\big)\big(G^{-1}(\bar{z}(s), \bar{y}(s))\big)\,ds +
      \sqrt{2\beta^{-1}} \big(\nabla
      \xi\,\sigma\big)_{il}\big(G^{-1}(\bar{z}(s), \bar{y}(s))\big)\, dw_l(s)\,,\\
      d\bar{y}_j(s) =& \big(\mathcal{L}\phi_j\big)\big(G^{-1}(\bar{z}(s),
      \bar{y}(s))\big)\,ds + \sqrt{2\beta^{-1}} \big(\nabla
      \phi\,\sigma\big)_{jl}\big(G^{-1}(\bar{z}(s), \bar{y}(s))\big)\, dw_l(s) \,,
    \end{split}
    \label{xi-phi-dynamics}
  \end{align}
  for $1 \le i \le m$ and $1 \le j \le n-m$, which can be considered as the
  equation of the dynamics $x(s)$ under the new coordinate $(\xi, \phi)$.
  Furthermore, the condition~(\ref{assump-xi-phi}) implies that the noise terms driving 
  the dynamics $\bar{z}(s)$ and $\bar{y}(s)$ in (\ref{xi-phi-dynamics}) are independent of each other.

The following result concerns the local existence of the function $\phi$. 
\vspace{0.1cm}

\begin{prop}
  Suppose Assumption~\ref{assump-0} holds. The matrix $a$ is $C^2$ smooth 
  and satisfies the condition (\ref{a-elliptic}) for some constant $c_1 >
  0$. The following two statements are equivalent. 
 \begin{enumerate}
   \item[(1)]
  There exists a
  neighborhood $\Omega_x$ of $x$, such that the map $G$ in (\ref{map-g}) is
  one to one from $\Omega_x$ to $\mbox{Im}(G)$, and that the
  condition~(\ref{assump-xi-phi}) is satisfied. 
\item[(2)]
There exists a neighborhood $\Omega_x$ of $x$, such that $\Pi^T B_{ij} \equiv 0$ on $\Omega_x$, for $1 \le i, j \le m$, where $\Pi$ is defined in (\ref{mat-pi}) and 
  $B_{ij} \in \mathbb{R}^n$ is given by 
  \begin{align}
  B_{ij, l'} = (a\nabla\xi_i)_l \frac{\partial(a\nabla\xi_j)_{l'}}{\partial x_l}
    - (a\nabla\xi_j)_l\frac{\partial(a\nabla\xi_i)_{l'}}{\partial x_l}\,,
    \quad 1 \le l' \le n\,.
\end{align}
  \end{enumerate}
  \label{prop-internal-coordinate-exist}
\end{prop}
\begin{proof}
  Let $u : \mathbb{R}^n \rightarrow \mathbb{R}$ be a $C^2$ smooth function.
  For $1 \le i \le m$, we define the differential operator $L_i$ by 
  \begin{align}
    L_i u = (a\nabla \xi_i)_l  \frac{\partial u}{\partial x_{l}}\,.
    \label{def-li}
  \end{align}
  By inverse mapping theorem, it is sufficient to find functions
  $\phi_1, \phi_2,\cdots, \phi_{n-m}$, which solve the PDE system 
  \begin{align}
    L_i \phi_j = 0 \,, \quad \mbox{for}~1 \le i \le m\,,
    \label{pde-li-u}
  \end{align}
  where $1\le j \le n-m$, such that $\nabla\phi_1, \nabla\phi_2, \cdots, \nabla\phi_{n-m}$ are linearly independent.
  From Frobenius theorem~\cite{lang1995differential}, such linearly independent
  solutions of the PDE system (\ref{pde-li-u}) exist if and only if there are
  functions $c_{ij}^k : \mathbb{R}^n \rightarrow \mathbb{R}$, such that 
  \begin{align}
  L_iL_j u - L_jL_i u = c_{ij}^k L_k u = c_{ij}^k (a\nabla \xi_k)_{l'} 
\frac{\partial u}{\partial x_{l'}}\,,
    \label{integral-condition}
  \end{align}
  holds for all $1 \le i, j \le m$, and for any $C^2$ function $u$.
  From (\ref{def-li}), we can directly compute 
  \begin{align}
    L_iL_j u - L_jL_i u = \bigg[(a\nabla\xi_i)_l \frac{\partial(a\nabla\xi_j)_{l'}}{\partial x_l}
    - (a\nabla\xi_j)_l\frac{\partial(a\nabla\xi_i)_{l'}}{\partial x_l}\bigg]
    \frac{\partial u}{\partial x_{l'}}
    =  B_{ij, l'} \frac{\partial u}{\partial x_{l'}}\,,
    \label{lij-lji-direct}
  \end{align}
    for $1 \le i, j \le m$. Now we prove the equivalence of the statements (1) and (2).

    $(1) \Rightarrow (2)$. 
    Suppose (\ref{integral-condition}) holds for some functions $c_{ij}^k$, then from
    (\ref{lij-lji-direct}) we have $B_{ij} = c_{ij}^k
    (a\nabla \xi_k)$. Since the matrix $\Pi$ satisfies $\Pi^Ta=a\Pi$ and
    $\Pi\nabla\xi_k=0$ in (\ref{mat-pi}), we conclude that 
    $$\Pi^T B_{ij} = c_{ij}^k (\Pi^Ta)\nabla\xi_k = c_{ij}^k a\Pi\nabla\xi_k =
    0\,.$$ 

    $(2) \Rightarrow (1)$. Suppose $\Pi^TB_{ij} \equiv 0$.  Then from the
    definition of $\Pi$ in (\ref{proj-pi}), we have 
    \begin{align}
      B_{ij} = (\Phi^{-1})_{kk'} \frac{\partial \xi_{k'}}{\partial x_{l'}} B_{ij,l'}
      (a\nabla \xi_k)\,, \quad 1 \le i, j \le m\,, 
    \end{align}
    which implies that (\ref{integral-condition}) holds if we choose $c_{ij}^k =
    (\Phi^{-1})_{kk'} \frac{\partial \xi_{k'}}{\partial x_{l'}} B_{ij,l'}$.
    Therefore the statement (1) is true by Frobenius theorem.
\end{proof}
\vspace{0.2cm}

\begin{remark}
  Proposition~\ref{prop-internal-coordinate-exist} provides conditions under
  which 
  we can reduce the case of a nonlinear reaction coordinate to 
  the linear reaction coordinate case in (\ref{xi-phi-dynamics}), locally in a
  neighborhood of a given state. The latter has been extensively investigated in
  literature in the study of slow-fast stochastic dynamical
  systems~\cite{pavliotis2008_multiscale,liu2010,asymptotic_analysis,khasminskii}.
  Although it seems impossible to solve $\phi$ for a general $\xi$ and matrix $a$
  provided that it exists,  
  it is interesting to mention the following special cases when $\phi$ exists
  or can be explicitly constructed.
  \begin{enumerate}
    \item
      When the reaction coordinate $\xi$ is scalar ($m=1$), 
      the statements of Proposition~\ref{prop-internal-coordinate-exist} are
      always true, i.e., the function $\phi$ always exists in this case.
    \item
      Consider $\xi(x) = (x_1, x_2, \cdots, x_m)^T$ is linear and the matrix $a =
      \mbox{diag}\{\sigma_1\sigma_1^T, \sigma_2\sigma_2^T\}$ is block diagonal, where 
      $\sigma_1\sigma_1^T \in \mathbb{R}^{m\times m}$ 
      and $\sigma_2\sigma_2^T \in \mathbb{R}^{(n-m)\times (n-m)}$.
      In this case, we can simply choose $\phi(x) = (x_{m+1}, \cdots, x_{n})^T$.
    \item
      Let $x=(x_1, x_2)^T $ be the state of a particle in
      $\mathbb{R}^2$ and $(r, \theta)$ denotes the
      coordinate of $x$ in the polar coordinate system. Assuming $a=I_{2\times
      2}$ and $\xi(x) = r = (x_1^2 + x_2^2)^{\frac{1}{2}}$, 
      we can verify that condition~(\ref{assump-xi-phi}) is satisfied with the function $\phi(x) =
      \theta(x)$. Furthermore, we note that this example can be generalized to the case of multiple
      particles where $\xi$ consists of radius or angles of different particles.
  \end{enumerate}
  \label{rmk-coordinate-transformation}
\end{remark}
\vspace{0.2cm}
 
In the following, let us assume that $\phi$ exists globally such
that the map $G$ in (\ref{map-g}) is one to one from $\mathbb{R}^n$ to itself. 
Given $z \in \mathbb{R}^m$, we consider the dynamics 
  \begin{align}
    \bar{x}(s) = G^{-1}(z, \bar{y}(s))\,, 
    \label{y-fix-z}
  \end{align} 
  where $\bar{y}(s)$ satisfies the second equation in (\ref{xi-phi-dynamics})
  with $\bar{z}(s)=z$ fixed. The following result states that the invariant measure of
    (\ref{y-fix-z}) coincides with $\mu_z$. 
\vspace{0.1cm}

 \begin{prop}
   Given $z \in \mathbb{R}^m$, the dynamics $\bar{x}(s)$ in (\ref{y-fix-z}) satisfies the SDE
  \begin{align}
    d\bar{x}_i(s) = -(\Pi^Ta)_{ij} \frac{\partial V}{\partial x_j}\, ds + \frac{1}{\beta}
    \frac{\partial (\Pi^Ta)_{ij}}{\partial x_j}\, ds + \sqrt{2\beta^{-1}}
    (\Pi^T\sigma)_{ij}\, dw_j(s)\,, \quad 1 \le i \le n\,.
  \end{align}
  In particular, $\bar{x}(s)\in \Sigma_z$ for $s \ge 0$ and it has a unique invariant
   measure $\mu_z$, which is defined in (\ref{measure-mu-z}). 
   \label{prop-phi-sde}
 \end{prop}
 \begin{proof}
   Clearly, (\ref{y-fix-z}) implies $\bar{x}(s) \in \Sigma_z$,
   for $\forall\,s \ge 0$. Applying Ito's formula to (\ref{y-fix-z}), we get 
  \begin{align}
    d\bar{x}_i(s) = \frac{\partial (G^{-1})_i}{\partial \phi_j} \mathcal{L}\phi_j\,ds +
    \frac{1}{\beta} (\nabla \phi a \nabla\phi^T)_{jl} \frac{\partial^2
    (G^{-1})_{i}}{\partial \phi_j\partial \phi_l} ds + \sqrt{2\beta^{-1}}
    \frac{\partial (G^{-1})_i}{\partial \phi_j} (\nabla \phi\,\sigma)_{jl}\, dw_l(s)\,,
    \label{y-fixed-xi}
  \end{align}
  where derivatives of $G^{-1}$ are evaluated at $\big(z, \phi(\bar{x}(s))\big)$,
  while functions $\mathcal{L}\phi_j$, $\nabla\phi a\nabla\phi^T$,
  and $\nabla\phi\,\sigma$ are evaluated at $\bar{x}(s)$. 

   Based on the discussions in Subsection~\ref{subsec-notations}, we
   know that, in order to prove the conclusion, it suffices to show the infinitesimal generator of
   (\ref{y-fixed-xi}) coincides with the operator $\mathcal{L}_0$ which is defined in (\ref{l-0-1}). For this purpose, taking derivatives in the identity
  \begin{align}
    G^{-1}(\xi(x), \phi(x)) = x\,, \quad \forall\,x \in \mathbb{R}^n\,,
  \end{align} 
    we have
  \begin{align}
    \frac{\partial (G^{-1})_i}{\partial \xi_l} \frac{\partial \xi_l}{\partial x_j} + 
    \frac{\partial (G^{-1})_i}{\partial \phi_l} \frac{\partial \phi_l}{\partial
    x_j}  = \delta_{ij} \,, \qquad \forall~1 \le i, j \le n\,.
    \label{chain-rule}
  \end{align}
  Together with (\ref{proj-pi}) and the condition~(\ref{assump-xi-phi}), we can obtain 
  \begin{align}
    \frac{\partial (G^{-1})_i}{\partial \xi_l} = (\Phi^{-1})_{ll'} (a\nabla \xi_{l'})_i\,,
  \end{align}
  as well as 
  \begin{align}
    \frac{\partial (G^{-1})_i}{\partial \phi_l} \frac{\partial \phi_l}{\partial
    x_j} = \delta_{ij} - \frac{\partial (G^{-1})_i}{\partial \xi_l}
    \frac{\partial \xi_l}{\partial
    x_j} = \delta_{ij} - 
    (\Phi^{-1})_{ll'} (a\nabla \xi_{l'})_i \frac{\partial \xi_l}{\partial x_j}
    = \Pi_{ji}\,.
  \end{align}
  For the first term on the right hand side of (\ref{y-fixed-xi}), using the expression
   (\ref{generator-l}) of $\mathcal{L}$ and noticing that the first argument of $G^{-1}$ is fixed, we can compute 
  \begin{align}
    \begin{split}
      \frac{\partial (G^{-1})_i}{\partial \phi_j} \mathcal{L}\phi_j
    =& 
\frac{e^{\beta V}}{\beta}
\frac{\partial (G^{-1})_i}{\partial \phi_j} 
    \frac{\partial}{\partial x_{i'}}\Big(a_{i'j'}
    e^{-\beta V}\frac{\partial \phi_j}{\partial x_{j'}}\Big) \\
    =&
\frac{e^{\beta V}}{\beta}
    \frac{\partial}{\partial x_{i'}}\bigg[
\frac{\partial (G^{-1})_i}{\partial \phi_j} 
    a_{i'j'}
    e^{-\beta V}\frac{\partial \phi_j}{\partial x_{j'}}\bigg] - \frac{1}{\beta}
  a_{i'j'} 
\frac{\partial^2 (G^{-1})_i}{\partial \phi_j\partial \phi_l} 
    \frac{\partial \phi_j}{\partial x_{j'}}
    \frac{\partial \phi_l}{\partial x_{i'}}\\
    =&
\frac{e^{\beta V}}{\beta}
    \frac{\partial}{\partial x_{i'}}\Big[(\Pi^Ta)_{ii'} e^{-\beta V}\Big] - \frac{1}{\beta}
  a_{i'j'} 
\frac{\partial^2 (G^{-1})_i}{\partial \phi_j\partial \phi_l} 
    \frac{\partial \phi_j}{\partial x_{j'}}
    \frac{\partial \phi_l}{\partial x_{i'}}\,.
    \end{split}
  \end{align}
  With the above computation, we know that the infinitesimal
  generator of $\bar{x}(s)$ in (\ref{y-fixed-xi}) is indeed $\mathcal{L}_0$,
  and the SDE (\ref{y-fixed-xi}) can be simplified as 
  \begin{align}
    d\bar{x}_i(s) = -(\Pi^Ta)_{ij} \frac{\partial V}{\partial x_j}\, ds + \frac{1}{\beta}
    \frac{\partial (\Pi^Ta)_{ij}}{\partial x_j}\, ds + \sqrt{2\beta^{-1}}
    (\Pi^T\sigma)_{ij}\, dw_j(s)\,.
  \end{align}
  Applying the result of~\cite[Theorem~4]{zhang2017}, we conclude that the invariant measure
  of the dynamics $\bar{x}(s)$ is given by $\mu_z$. 
\end{proof}

  \section{Proof of Lemma~\ref{lemma-spectral-gap-stiff-v}}
  \label{sec-app-2}
  This appendix is devoted to proving Lemma~\ref{lemma-spectral-gap-stiff-v}.
  We will only sketch the proof, since we essentially follow the argument in \cite{bakry2013analysis} 
(also see \cite[Chap.  14]{villani2008optimal}) with some technical modifications. Before entering
  the proof, we need to first introduce some notations. 

  In the following, for fixed $z \in \mathbb{R}^m$, we will denote by $N$ the Riemannian submanifold $\Sigma_z$
  where the metric is induced from the Euclidean distance on $\mathbb{R}^n$.  Let $\nabla^N$, $\Delta^N$ be the gradient operator and the Laplacian
  operator on $N$, respectively. Recalling the parameter $\epsilon \ll 1$ and the potential
  function $V_1$ in (\ref{v-v0-v1}), we consider the operator  
  \begin{align}
    \mathcal{L}^N = -\frac{1}{\epsilon} \nabla^NV_1 \cdot \nabla^N + \frac{1}{\beta}
    \Delta^N
  \end{align}
  on $N$ and denote by $(T_t)_{t\ge 0}$ the corresponding semigroup. It is straightforward to
  verify that $T_t$ is invariant with respect to the probability measure
  $\bar{\nu}$ which is given by
  \begin{align}
    d\bar{\nu} = \frac{1}{Z} e^{-\frac{\beta}{\epsilon} V_1} d\nu_z\,,
    \label{bar-nu}
  \end{align}
  where $Z$ is the normalization constant and $\nu_z$ denotes the surface measure on $N$.

Given two smooth functions $f,\, h : N \rightarrow \mathbb{R}$,
the associated $\Gamma$ operator (\textit{carr\'{e} du champ}) and $\Gamma_2$
operator of $\mathcal{L}^N$ are defined as 
  \begin{align}
    \begin{split}
    \Gamma(f,h) =& \frac{1}{2} \Big[\mathcal{L}^N(fh) - f\mathcal{L}^Nh - h\mathcal{L}^Nf\Big]
    = \frac{1}{\beta}\, \nabla^Nf\cdot\nabla^Nh \,,\\
    \Gamma_2(f,h) =& \frac{1}{2} \Big[\mathcal{L}^N\Gamma(f,h) -
    \Gamma(f,\mathcal{L}^Nh) - \Gamma(\mathcal{L}^Nf, h)\Big]\,.
  \end{split}
  \label{gamma-1-2}
  \end{align}
   Let us consider the (smooth) extensions of $f,h$ from $N$ to
   $\mathbb{R}^n$, which we again denote by $f$ and $h$, respectively. Recall that $P$ is the orthogonal projection
   operator from $T_x\mathbb{R}^n$ to $T_xN$ introduced in Subsection~\ref{subsec-notations}. 
   We can check that $\nabla^N f= P\nabla f$, $\nabla^N h = P\nabla h$, and
  therefore from (\ref{gamma-1-2}) we have
  \begin{align}
    \Gamma(f,h) = \frac{1}{\beta} (P\nabla f)\cdot (P \nabla h)\,.
    \label{gamma-f-h}
  \end{align}
  Clearly, the above
  expression of $\Gamma$ does not depend on the extensions of $f$ and $h$ we choose. 

 For the $\Gamma_2$ operator in (\ref{gamma-1-2}), applying the Bochner-Lichnerowicz formula~\cite[Theorem C.3.3]{bakry2013analysis},
  we can compute 
  \begin{align}
    \begin{split}
    \Gamma_2(f,f) =\,& 
     \frac{1}{2} \mathcal{L}^N\Gamma(f,f) - \Gamma(f, \mathcal{L}^Nf) \\
     =\,& \frac{1}{\beta^2}\|\mbox{Hess}^N f\|^2_{HS} + \Big(\frac{1}{\beta^2}\mbox{Ric}^N +
     \frac{1}{\epsilon\beta}\mbox{Hess}^N V_1\Big) (\nabla^N f,\nabla^N f)\\
     \ge\, & \Big(\frac{1}{\beta^2}\mbox{Ric}^N +
     \frac{1}{\epsilon\beta}\mbox{Hess}^N V_1\Big) (\nabla^N f,\nabla^N f)\,.
  \end{split}
  \label{gamma-2-bound-ric-hessian}
  \end{align}
  In the above, $\|\mbox{Hess}^Nf\|_{HS}$ is the Hilbert-Schmidt norm of the Hessian of
  the function $f$, and $\mbox{Ric}^N$ denotes the Ricci tensor on $N$.

  After the above preparations, we are ready to prove Lemma~\ref{lemma-spectral-gap-stiff-v}.
  
  \vspace{0.1cm}
  \begin{proof}[Proof of Lemma~\ref{lemma-spectral-gap-stiff-v}]
  We divide the proof into two steps.
  \begin{enumerate}
    \item 
Firstly, let us prove the Poincar{\'e} inequality for the invariant measure $\bar{\nu}$, i.e., 
  \begin{align}
    \int_N f^2\,d\bar{\nu} - \Big(\int_N f\,d\bar{\nu}\Big)^2 \le
    \frac{2\epsilon}{K} \int_N \Gamma(f,f)\, d\bar{\nu}\,,
    \label{poincare-nu-bar}
  \end{align}
for all smooth functions $f : N \rightarrow \mathbb{R}$, when $\epsilon$ is small enough.
       According to~\cite[Proposition 4.8.1]{bakry2013analysis}, 
      it is sufficient to prove the curvature condition
      $CD\big(\frac{K}{2\epsilon}, \infty\big)$, i.e.,  
      \begin{align}
  \Gamma_2(f,f) \ge \frac{K}{2\epsilon} \Gamma(f,f)\,,
  \label{cd-condition}
\end{align}
for all smooth functions $f : N \rightarrow \mathbb{R}$.
      Notice that, the $K$-\,convexity and $C^2$ smoothness of $V_1$ imply
  \begin{align}
  \mbox{Hess}^NV_1(\nabla^Nf,\nabla^Nf) \ge K\,|\nabla^Nf|^2\,.
  \label{hess-bound}
\end{align}
      Denote by $R^N$, $H$ the Riemannian curvature tensor and 
      the mean curvature vector of $N$, respectively. Given $x \in N$, let $\bm{e}_i \in T_xN$, $1 \le i \le n-m$, be an
      orthonormal basis of $T_xN$.  
      Applying the Gauss equation \cite[Theorem~8.4]{lee1997riemannian} and using
      the relation $H=(I-P)\sum_{i=1}^{n-m}\nabla_{\bm{e}_i}\bm{e}_i$ \cite[Proposition~1]{zhang2017}, we can compute
\begin{align}
  \begin{split}
  &\mbox{Ric}^N(X, X) \\
  =& \sum_{i=1}^{n-m} \Big[R^N(X, \bm{e}_i)\bm{e}_i\Big] \cdot X \\
  =& \sum_{i=1}^{n-m}\Big[- \Big((I-P)\nabla_X\bm{e}_i\Big) \cdot \Big((I-P) \nabla_{\bm{e_i}} X\Big)
  + \Big((I-P) \nabla_{\bm{e}_i}\bm{e}_i\Big) \cdot\Big((I-P) \nabla_{X}X\Big)\Big]\\
  =& -\sum_{i=1}^{n-m}\Big((I-P)\nabla_X\bm{e}_i\Big)\cdot\Big((I-P) \nabla_{X}\bm{e_i}\Big)
  - (\nabla_XH) \cdot X\, ,
\end{split}
  \label{ric-expressions}
\end{align}
for all $X \in T_xN$. 
      In the above, we have used the fact that
      the Riemannian curvature tensor of the Euclidean space $\mathbb{R}^n$ vanishes, as
      well as 
      \begin{align*}
      (I-P)\big(\nabla_{\bm{e}_i} X- \nabla_X\bm{e}_i\big) = (I-P)[\bm{e}_i, X] = 0\,,
      \end{align*}
      and $X\cdot H = 0$, since $\bm{e}_i, X \in T_xN$ and $H \in (T_xN)^{\perp}$. 
      From the last expression in (\ref{ric-expressions}) and the assumptions in Lemma~\ref{lemma-spectral-gap-stiff-v}, it is not difficult to conclude that $\exists C\in
\mathbb{R}$, such that 
\begin{align}
  \mbox{Ric}^N(X,X) \ge C|X|^2\,,  \quad \forall~ X \in T_xN\,.
  \label{ric-bound}
\end{align}
Combining (\ref{gamma-2-bound-ric-hessian}), (\ref{hess-bound}), and
(\ref{ric-bound}), we obtain 
\begin{align}
  \Gamma_2(f,f) \ge \frac{1}{\beta} \Big(\frac{C}{\beta} +
  \frac{K}{\epsilon}\Big) |\nabla^N f|^2 \ge \frac{K}{2\epsilon}
  \Gamma(f,f)\,,
  \label{gamma-2-bound-by-gamma-1}
\end{align}
      when $\epsilon$ is small enough. Therefore, the curvature condition
      (\ref{cd-condition}) is satisfied and the Poincar{\'e} inequality
      (\ref{poincare-nu-bar}) follows.
\item
    Secondly, we derive the Poincar{\'e} inequality for the measure $\mu_z$ using
  Holley-Stroock perturbation lemma~\cite{Holley1987,Ledoux2001}. 
   For this purpose, from (\ref{measure-mu-z}) and (\ref{bar-nu}), we know that the probability
   measure $\mu_z$ is related to $\bar{\nu}$ by
  \begin{align}
    d\mu_z = \frac{1}{Z} e^{-\beta V_0} \Big[\mbox{det}(\nabla\xi \nabla\xi^T)\Big]^{-\frac{1}{2}} d\bar{\nu}\,,
  \end{align}
  where $Z$ is the normalization constant.
  And our assumptions imply that both $\frac{d\mu_z}{d\bar{\nu}}$ and
  $\frac{d\bar{\nu}}{d\mu_z}$ are bounded on $N$ by some constant $C>0$.
  Therefore,
  applying~\cite[Proposition~4.2.7]{bakry2013analysis}, we have 
  \begin{align}
    \begin{split}
     \int_{\Sigma_z} f^2\, d\mu_z - \Big(\int_{\Sigma_z} f\, d\mu_z\Big)^2 
    \le \frac{2\epsilon C}{K} \int_{\Sigma_z} \Gamma(f,f)\, d\mu_z 
      = \frac{2\epsilon C}{\beta K} \int_{\Sigma_z} \, |P\nabla f|^2 d\mu_z \,,
    \end{split}
    \label{nu-to-mu}
  \end{align} 
 where the constant $C$ may differ from the upper bound of
 $\frac{d\bar{\nu}}{d\mu_z}$ and $\frac{d\mu_z}{d\bar{\nu}}$.
Assuming that $f$ has been extended from $N$ to $\mathbb{R}^n$, 
(\ref{a-elliptic}) and (\ref{mat-pi}) imply 
\begin{align}
  |P\nabla f|^2 \le |\Pi P\nabla f|^2 =
  |\Pi \nabla f|^2 \le \frac{1}{c_1} (a\Pi \nabla f)\cdot
  (\Pi\nabla f)\,.
\end{align}
Therefore, the inequality (\ref{poincare-ineq-prop}) follows readily from (\ref{nu-to-mu}) and the expression of the
Dirichlet form $\mathcal{E}_z$ in (\ref{dirichlet-form-z}).
  \end{enumerate}
\end{proof}
  \section{Mean square error estimate of marginals}
  \label{sec-app-3}
  In this appendix, instead of assuming the Lipschitz condition on $\widetilde{b}$ (Assumption~\ref{assump-1}), 
we provide a mean square error estimate of the marginals for the effective
dynamics and the process $\xi(x(s))$, under the following dissipative assumption. 
      \vspace{0.1cm}

      \begin{assump}
	$\exists$ $L_d, L_\sigma > 0$, such that $\forall\, z, z' \in
	\mathbb{R}^m$, we have 
      \begin{align}
	 \big(\widetilde{b}(z) - \widetilde{b}(z')\big)\cdot(z - z') \le - L_d |z
	 - z'|^2\,,\qquad 
	 \big\|\widetilde{\sigma}(z) - \widetilde{\sigma}(z')\big\|_F \le
    L_\sigma|z-z'|\,. 
    \label{b-dissipative-sigma-lip}
  \end{align}
    \label{assump-dissipative}
      \end{assump}
      \begin{prop}
	Suppose that Assumptions~\ref{assump-0},\,\ref{assump-2},\,\ref{assump-4} and~\ref{assump-dissipative} hold.
$x(s)$ satisfies the SDE (\ref{dynamics-a}) starting from
  $x(0)\sim \mu$, and $z(s)$ is the effective dynamics (\ref{modified-eff-dynamics}) with $z(0) =
  \xi(x(0))$.
	\begin{enumerate}
	  \item
Assume $L_d > \frac{L_\sigma^2}{\beta}$. 
      Choose $v_1, v_2>0$ such that $$C_1= L_d -
	    \frac{L_\sigma^2(1+v_2)}{\beta}-\frac{v_1}{2} > 0\,.$$ 
We have 
      \begin{align}
	\mathbf{E}\big|\xi(x(t)) - z(t)\big|^2 \le
	\frac{C_1^{-1}}{\beta\rho}\bigg[\frac{\kappa_1^2}{2v_1} +
	\frac{2\kappa_2^2}{\beta}\Big(1+\frac{1}{v_2}\Big)\bigg]\Big(1-e^{-2C_1t}\Big) \,, \quad
	\forall\, t\ge 0\,.
      \end{align}
	  \item
Assume $L_d \le \frac{L_\sigma^2}{\beta}$. For any $v_1, v_2>0$, we define
$$C_2=  \frac{L_\sigma^2(1+v_2)}{\beta}-L_d + \frac{v_1}{2} > 0\,.$$ 
We have 
      \begin{align}
	\mathbf{E}\big|\xi(x(t)) - z(t)\big|^2 \le
	\frac{C_2^{-1}}{\beta\rho}\bigg[\frac{\kappa_1^2}{2v_1} +
	\frac{2\kappa_2^2}{\beta}\Big(1+\frac{1}{v_2}\Big)\bigg]\Big(e^{2C_2t}-1\Big) \,, \quad
	\forall\, t\ge 0\,.
      \end{align}
	\end{enumerate}
      \end{prop}
      \begin{proof} 
Recall the function $\varphi$ defined in (\ref{f-lxi-b-diff}).
	Using (\ref{xi-z-diff}) and applying Ito's formula, we obtain
\begin{align}
  \begin{split}
    \frac{1}{2}\big|\xi(x(t)) - z(t)\big|^2 =& \int_0^t
    \varphi(x(s))\cdot\big(\xi(x(s))-z(s)\big)\,ds + \int_0^t
    \Big(\widetilde{b}\big(\xi(x(s))\big) - \widetilde{b}(z(s))\Big)\cdot
    \big(\xi(x(s))-z(s)\big)\,ds \\
    &+ \frac{1}{\beta} \int_0^t \big\|A(x(s)) - \widetilde{\sigma}(z(s))\big\|_F^2\,ds
    + \sqrt{2\beta^{-1}} M(t)\,,
  \end{split}
  \label{xi-z-sqr-diff}
\end{align}
where  
\begin{align*}
  M(t) = \int_0^t \big(\xi(x(s)) - z(s)\big)^T \Big(A\big(x(s)\big) -
  \widetilde{\sigma}\big(z(s)\big)\Big)\,d\widetilde{w}_s
\end{align*}
 is the martingale term.
 Taking expectation in (\ref{xi-z-sqr-diff}) and differentiating with respect to
 time $t$, we get 
 \begin{align}
   \begin{split}
     \frac{1}{2}\frac{d}{dt}\mathbf{E}\big|\xi(x(t)) - z(t)\big|^2 
    = & 
     \mathbf{E}\Big[\varphi(x(t))\cdot \big(\xi(x(t))-z(t)\big)\Big]+
     \mathbf{E}\Big[\Big(\,\widetilde{b}\big(\xi(x(t))\big) -
     \widetilde{b}\big(z(t)\big)\Big)\cdot\big(\xi(x(t))-z(t)\big)\Big]\,\\
    &+ \frac{1}{\beta}\mathbf{E}\big\|A(x(t)) -
    \widetilde{\sigma}(z(t))\big\|_F^2\,. 
  \end{split}
  \label{app-marginal-proof-eq1}
 \end{align}
 Applying Assumption~\ref{assump-dissipative},
	Lemmas~\ref{lemma-f2-by-poincare}--\ref{lemma-Frobenius-estimate}, together with Young's
	inequality, we can estimate the right hand side of (\ref{app-marginal-proof-eq1}) and obtain
      \begin{align*}
  \begin{split}
    & \frac{1}{2}\frac{d}{dt}\mathbf{E}\big|\xi(x(t)) - z(t)\big|^2 \\
    \le & \frac{1}{2v_1} \mathbf{E}_\mu |\varphi|^2
    -\Big(L_d-\frac{v_1}{2}\Big) 
    \mathbf{E}\big|\xi(x(t)) - z(t)\big|^2 + \frac{2\kappa_2^2(1+\frac{1}{v_2})}{\beta^2
    \rho} + \frac{L_\sigma^2(1+v_2)}{\beta} \mathbf{E}\big|\xi(x(t)) - z(t)\big|^2\\
    \le & \frac{1}{\beta\rho}\bigg[\frac{\kappa_1^2}{2v_1} +
    \frac{2\kappa_2^2}{\beta}\Big(1+\frac{1}{v_2}\Big)\bigg] +
    \Big(\frac{L_\sigma^2(1+v_2)}{\beta}-L_d+\frac{v_1}{2}\Big)
    \mathbf{E}\big|\xi(x(t)) - z(t)\big|^2\,,
  \end{split}
\end{align*}
for any $v_1, v_2 > 0$. The conclusions follow by applying Gronwall's inequality. 
    \end{proof}
\bibliographystyle{siam}
\bibliography{reference}

\begin{thebibliography}{10}

\bibitem{ANDO1979}
{\sc T.~Ando}, {\em Concavity of certain maps on positive definite matrices and
  applications to {Hadamard} products}, Linear Algebra Appl., 26 (1979),
  pp.~203--241.

\bibitem{weighted-nash}
{\sc D.~Bakry, F.~Bolley, I.~Gentil, and P.~Maheux}, {\em Weighted {Nash}
  inequalities}, Rev. Mat. Iberoam., 28 (2012), pp.~879--906.

\bibitem{bakry2013analysis}
{\sc D.~Bakry, I.~Gentil, and M.~Ledoux}, {\em Analysis and Geometry of Markov
  Diffusion Operators}, Grundlehren der mathematischen Wissenschaften, Springer
  International Publishing, 2014.

\bibitem{banyaga2004lectures}
{\sc A.~Banyaga and D.~Hurtubise}, {\em Lectures on Morse Homology}, Texts in
  the Mathematical Sciences, Springer Netherlands, 2004.

\bibitem{asymptotic_analysis}
{\sc A.~Bensoussan, J.~L. Lions, and G.~Papanicolaou}, {\em Asymptotic analysis
  for periodic structures}, Studies in mathematics and its applications,
  North-Holland, 1978.

\bibitem{projection_diffusion}
{\sc G.~Ciccotti, T.~Leli{\`e}vre, and E.~Vanden-Eijnden}, {\em Projection of
  diffusions on submanifolds: Application to mean force computation}, Commun.
  Pur. Appl. Math., 61 (2008), pp.~371--408.

\bibitem{davies1990heat}
{\sc E.~B. Davies}, {\em Heat Kernels and Spectral Theory}, Cambridge Tracts in
  Mathematics, Cambridge University Press, 1990.

\bibitem{mimick_sde}
{\sc I.~Gy{\"o}ngy}, {\em Mimicking the one-dimensional marginal distributions
  of processes having an {Ito} differential}, Probab. Theory Related Fields.,
  71 (1986), pp.~501--516.

\bibitem{zhws16}
{\sc C.~Hartmann, Ch. Sch{\"u}tte, M.~Weber, and W.~Zhang}, {\em Importance
  sampling in path space for diffusion processes with slow-fast variables},
  Probab. Theory Related Fields., 170 (2018), pp.~177--228.

\bibitem{Holley1987}
{\sc R.~Holley and D.~Stroock}, {\em Logarithmic {Sobolev} inequalities and
  stochastic ising models}, J. Stat. Phys., 46 (1987), pp.~1159--1194.

\bibitem{peptide_cmd}
{\sc G.~Hummer and I.~G. Kevrekidis}, {\em Coarse molecular dynamics of a
  peptide fragment: Free energy, kinetics, and long-time dynamics
  computations}, J. Chem. Phys., 118 (2003), pp.~10762--10773.

\bibitem{kevrekidis2003}
{\sc I.~G. Kevrekidis, C.~W. Gear, J.~M. Hyman, P.~G Kevrekidid, O.~Runborg,
  and C.~Theodoropoulos}, {\em Equation-free, coarse-grained multiscale
  computation: Enabling mocroscopic simulators to perform system-level
  analysis}, Commun. Math. Sci., 1 (2003), pp.~715--762.

\bibitem{khasminskii}
{\sc R.~Khasminskii}, {\em Principle of averaging for parabolic and elliptic
  differential equations and for {M}arkov processes with small diffusion},
  Theory Probab. Appl., 8 (1963), pp.~1--21.

\bibitem{komorowski2012fluctuations}
{\sc T.~Komorowski, C.~Landim, and S.~Olla}, {\em Fluctuations in Markov
  Processes: Time Symmetry and Martingale Approximation}, Grundlehren der
  mathematischen Wissenschaften, Springer Berlin Heidelberg, 2012.

\bibitem{lang1995differential}
{\sc S.~Lang}, {\em Differential and Riemannian Manifolds}, Graduate Texts in
  Mathematics, Springer-Verlag, 1995.

\bibitem{Ledoux2001}
{\sc M.~Ledoux}, {\em Logarithmic {Sobolev} inequalities for unbounded spin
  systems revisited}, in S{\'e}minaire de Probabilit{\'e}s XXXV, J.~Az{\'e}ma,
  M.~{\'E}mery, M.~Ledoux, and M.~Yor, eds., Springer Berlin Heidelberg, 2001,
  pp.~167--194.

\bibitem{lee1997riemannian}
{\sc J.M. Lee}, {\em Riemannian Manifolds: An Introduction to Curvature},
  Graduate Texts in Mathematics, Springer New York, 1997.

\bibitem{effective_dynamics}
{\sc F.~Legoll and T.~Leli{\`e}vre}, {\em Effective dynamics using conditional
  expectations}, Nonlinearity, 23 (2010), pp.~2131--2163.

\bibitem{LEGOLL2017-pathwise}
{\sc F.~Legoll, T.~Leli{\`e}vre, and S.~Olla}, {\em Pathwise estimates for an
  effective dynamics}, Stoch. Process. Appl., 127 (2017), pp.~2841--2863.

\bibitem{LIEB1973}
{\sc E.~H. Lieb}, {\em Convex trace functions and the {Wigner-Yanase-Dyson}
  conjecture}, Adv. Math., 11 (1973), pp.~267--288.

\bibitem{liu2010}
{\sc D.~Liu}, {\em Strong convergence of principle of averaging for multiscale
  stochastic dynamical systems}, Commun. Math. Sci., 8 (2010), pp.~999--1020.

\bibitem{lyons1994}
{\sc T.~J. Lyons and T.~S. Zhang}, {\em Decomposition of {Dirichlet} processes
  and its application}, Ann. Probab., 22 (1994), pp.~494--524.

\bibitem{Maragliano2006}
{\sc L.~Maragliano and E.~Vanden-Eijnden}, {\em A temperature accelerated
  method for sampling free energy and determining reaction pathways in rare
  events simulations}, Chem. Phys. Lett., 426 (2006), pp.~168--175.

\bibitem{Mattingly2002}
{\sc J.~C. Mattingly, A.~M. Stuart, and D.~J. Higham}, {\em Ergodicity for
  {SDEs} and approximations: locally {Lipschitz} vector fields and degenerate
  noise}, Stoch. Proc. Appl., 101 (2002), pp.~185--232.

\bibitem{oksendalSDE}
{\sc B.~{\O}ksendal}, {\em Stochastic Differential Equations: An Introduction
  with Applications}, Springer, 5th~ed., 2000.

\bibitem{pavliotis2008_multiscale}
{\sc G.~A. Pavliotis and A.~M. Stuart}, {\em Multiscale Methods: Averaging and
  Homogenization}, Texts in Applied Mathematics, Springer New York, 2008.

\bibitem{Sharma2017}
{\sc U.~Sharma}, {\em Coarse-graining of {Fokker-Planck} equations}, PhD
  thesis, Technische Universiteit Eindhoven, 2017.

\bibitem{vanden-eijnden2003}
{\sc E.~Vanden-Eijnden}, {\em Numerical techniques for multi-scale dynamical
  systems with stochastic effects}, Commun. Math. Sci., 1 (2003), pp.~385--391.

\bibitem{revisit-ftstring}
{\sc E.~Vanden-Eijnden and M.~Venturoli}, {\em Revisiting the finite
  temperature string method for the calculation of reaction tubes and free
  energies}, J. Chem. Phys., 130 (2009), p.~194103.

\bibitem{villani2008optimal}
{\sc C.~Villani}, {\em Optimal Transport: Old and New}, Grundlehren der
  mathematischen Wissenschaften, Springer Berlin Heidelberg, 2008.

\bibitem{zhang2017}
{\sc W.~Zhang}, {\em Constructing ergodic diffusion processes on submanifolds},
  submitted,  (2017).

\bibitem{effective_dyn_2017}
{\sc W.~Zhang, C.~Hartmann, and Ch. Sch{\"u}tte}, {\em Effective dynamics along
  given reaction coordinates{,} and reaction rate theory}, Faraday Discuss.,
  195 (2016), pp.~365--394.

\end{thebibliography}
\end{document}